\documentclass[11pt]{article}
\usepackage{amsmath, amsthm, amssymb}    
\usepackage{bridges}                     
\usepackage{graphicx}                    
\usepackage{hyperref}                    

\usepackage{enumitem}
\usepackage{color}
\usepackage[style=base]{caption}
\usepackage{subfig}
\usepackage{wrapfig}
\usepackage{pdflscape}
\usepackage{array}
\usepackage{MnSymbol}

\newcommand{\reffig}[1]{Figure~\ref{Fig:#1}}

\newcommand{\cover}[1]{{\widetilde{#1}}}
\newcommand{\bdy}{\partial}

\newcommand{\calT}{\mathcal{T}}

\newcommand{\HH}{\mathbb{H}}
\newcommand{\ZZ}{\mathbb{Z}}

\title{Cohomology fractals}
\author{David Bachman\\
Pitzer College\\
{\tt bachman@pitzer.edu}
\and
Saul Schleimer\\
University of Warwick\\
{\tt s.schleimer@warwick.ac.uk}
\and
Henry Segerman\\
Oklahoma State University\\
{\tt henry@segerman.org}
}

\date{}				

\begin{document}

\maketitle

\thispagestyle{empty}

\begin{abstract}
We introduce cohomology fractals; these are certain images associated to a cohomology class on a hyperbolic three-manifold.  They include images made entirely from circles, and also images with no geometrically simple features.  They are closely related to limit sets of kleinian groups, but have some key differences.  As a consequence, we can zoom in almost any direction to arbitrary depth in real time.  We present an implementation in the setting of ideal triangulations using ray-casting.
\end{abstract}




\begin{figure}[htb!]
\centering
\subfloat[\texttt{m004}]{
\includegraphics[width=0.32\textwidth]{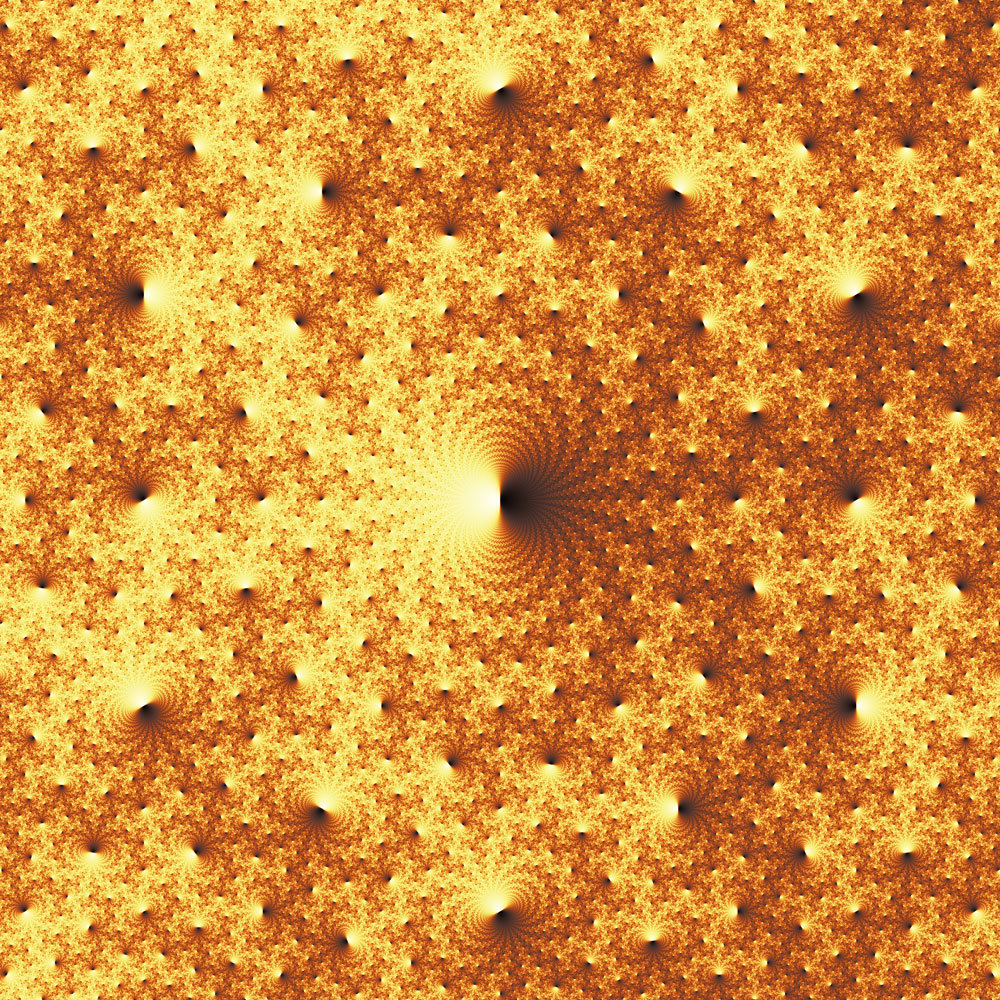}\label{Fig:m004}
}
\subfloat[\texttt{m129}]{
\includegraphics[width=0.32\textwidth]{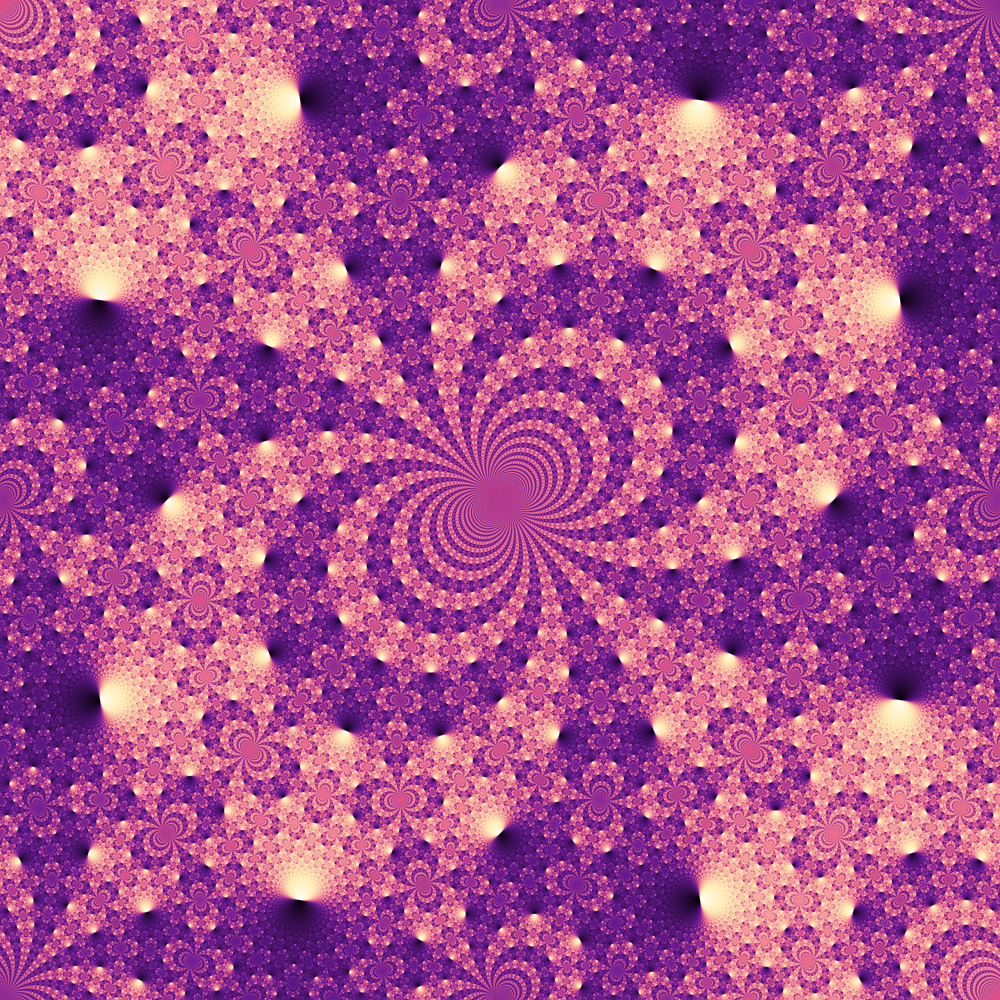}\label{Fig:m129}
}
\subfloat[\texttt{m345}]{
\includegraphics[width=0.32\textwidth]{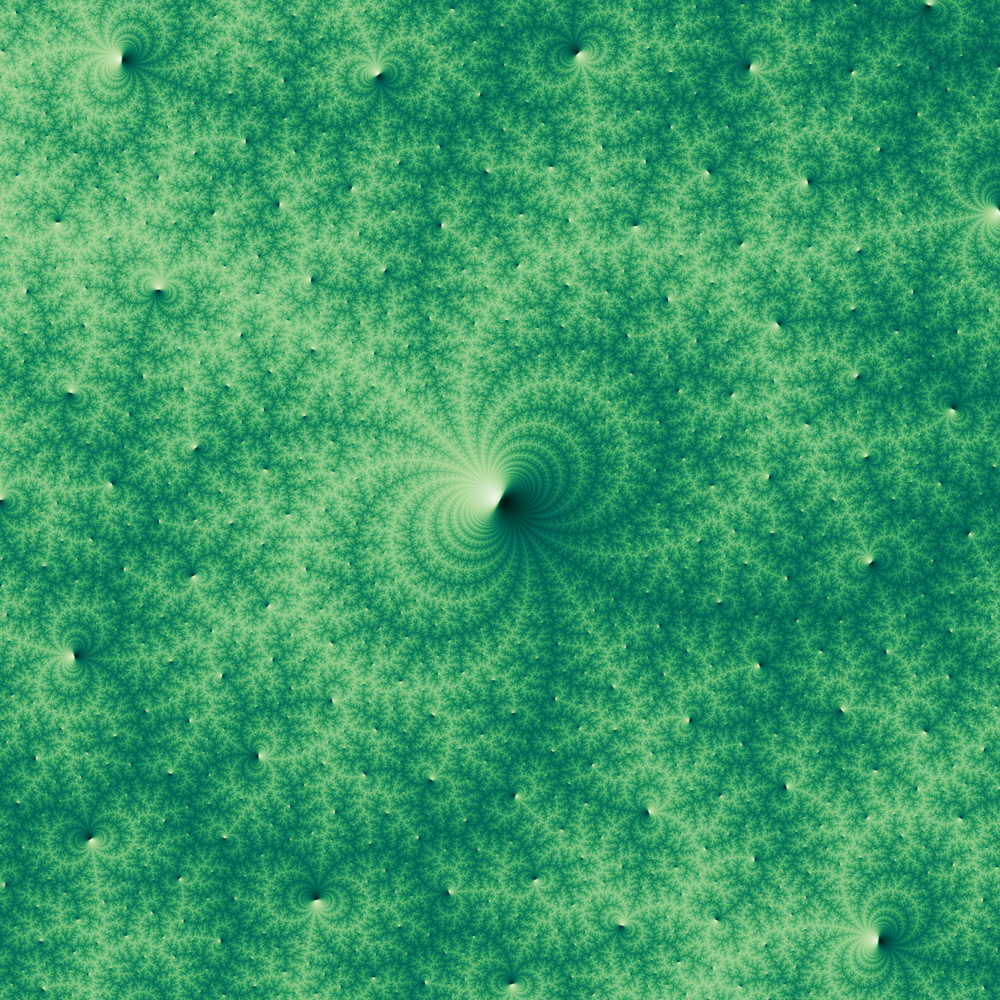}\label{Fig:m345}
}

\subfloat[\texttt{s000}]{
\includegraphics[width=0.32\textwidth]{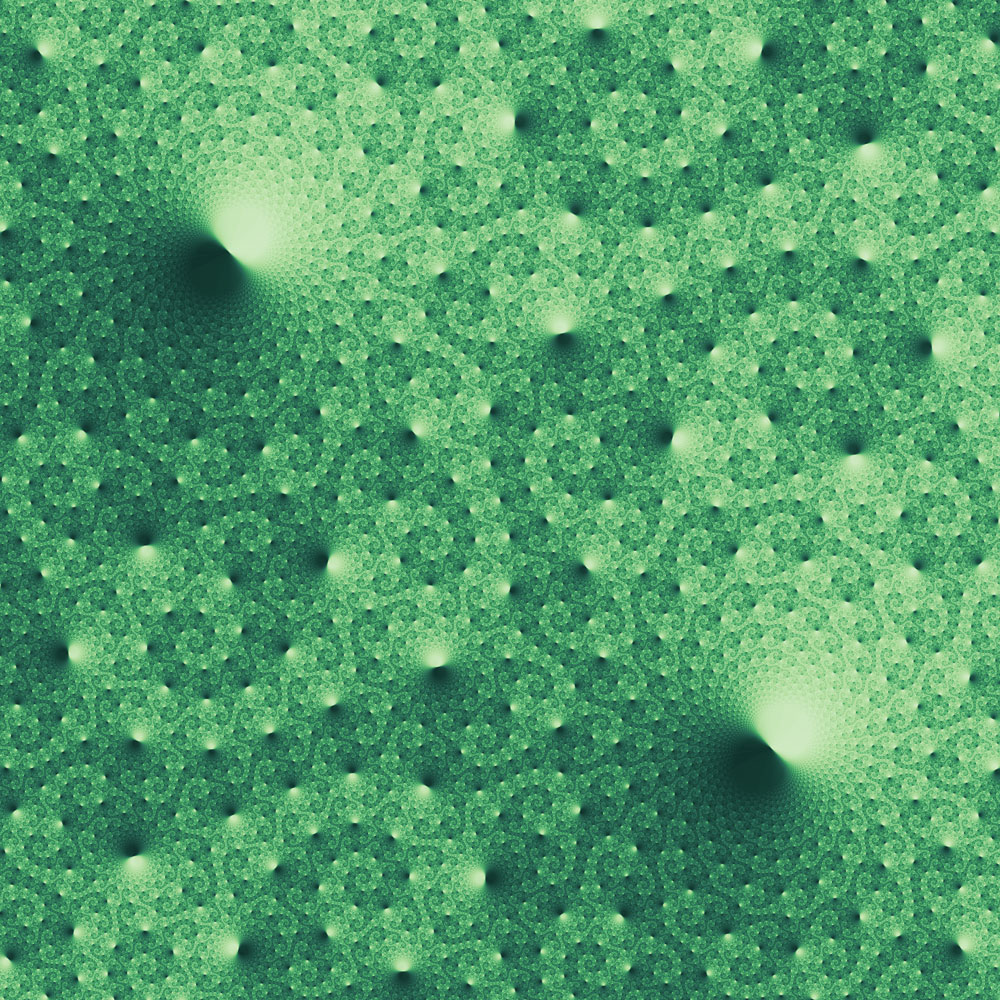}\label{Fig:s000}
}
\subfloat[\texttt{s227}]{
\includegraphics[width=0.32\textwidth]{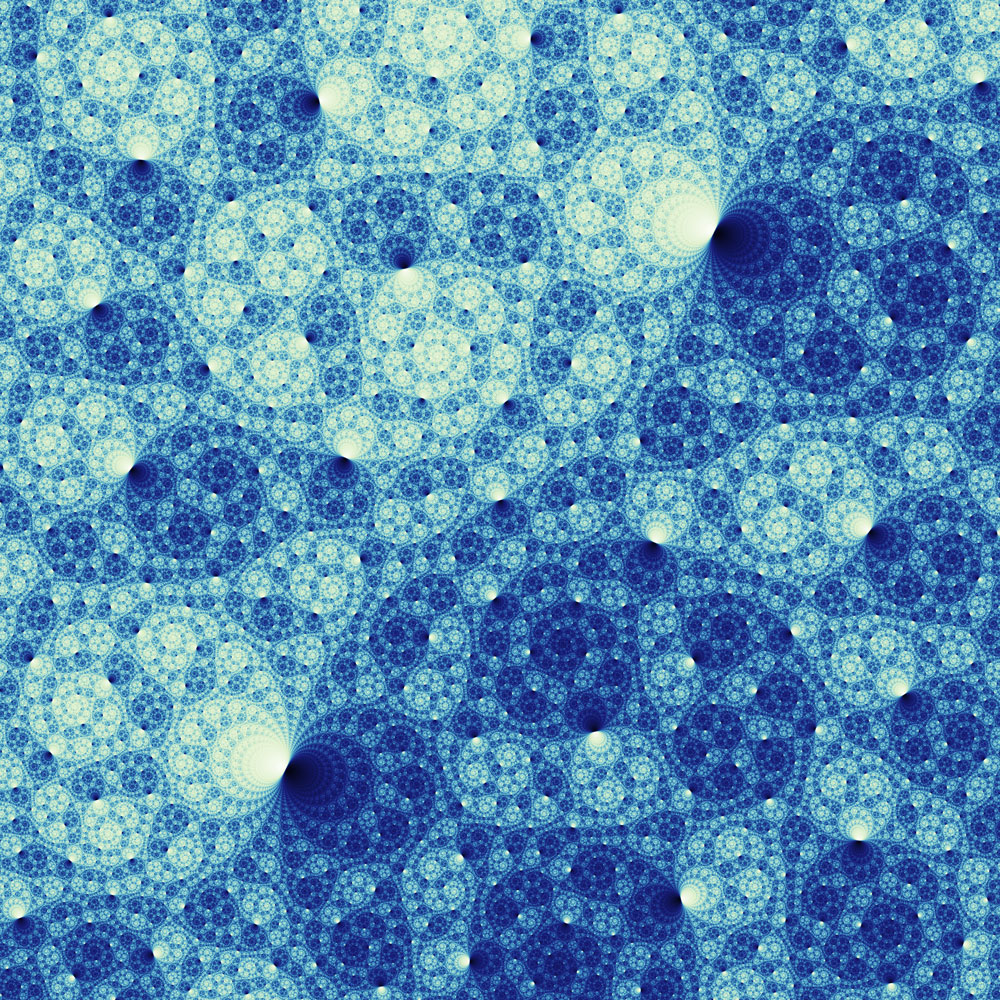}\label{Fig:s227}
}
\subfloat[\texttt{s776}]{
\includegraphics[width=0.32\textwidth]{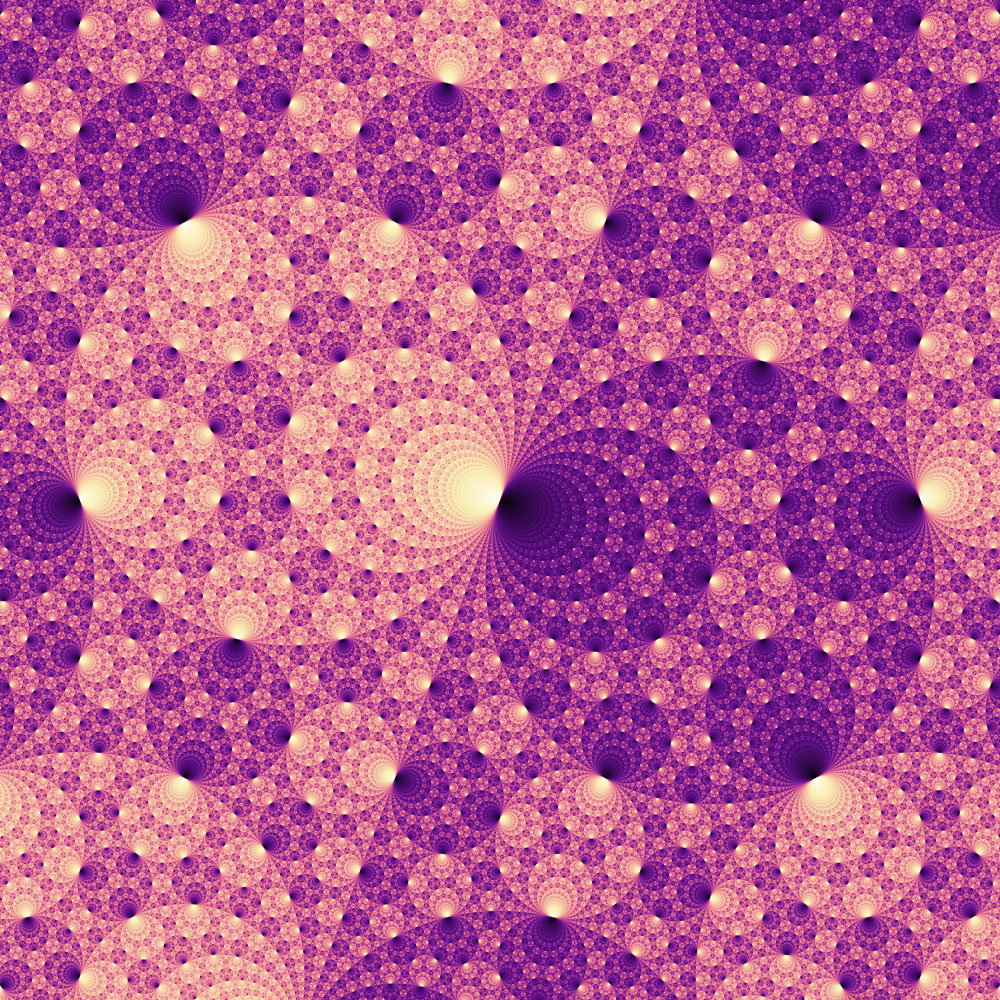}\label{Fig:s776}
}
\caption{\emph{Cohomology fractals. We list the names of the corresponding manifolds in the} \texttt{SnapPy}~\cite{snappy} \emph{census. The choice of colour scheme for each image is aesthetic: the mathematical content is the brightness of each pixel.}}
\label{Fig:PixelImages}
\end{figure}

We show a few examples of cohomology fractals in \reffig{PixelImages}.  We provide an online application that generates these, and many others, at \url{https://henryseg.github.io/cohomology_fractals/}.  The code is available at \url{https://github.com/henryseg/cohomology_fractals}.  We have also made YouTube videos explaining the basic ideas and giving a collection of zooms.  They are available at \url{https://youtu.be/fhBPhie1Tm0} and at \url{https://youtu.be/-g1wNbC9AxI}.

\begin{figure}[htb]
\centering
\subfloat[\emph{A torus.}]{
\label{Fig:Torus}
\includegraphics[height=1.67in]{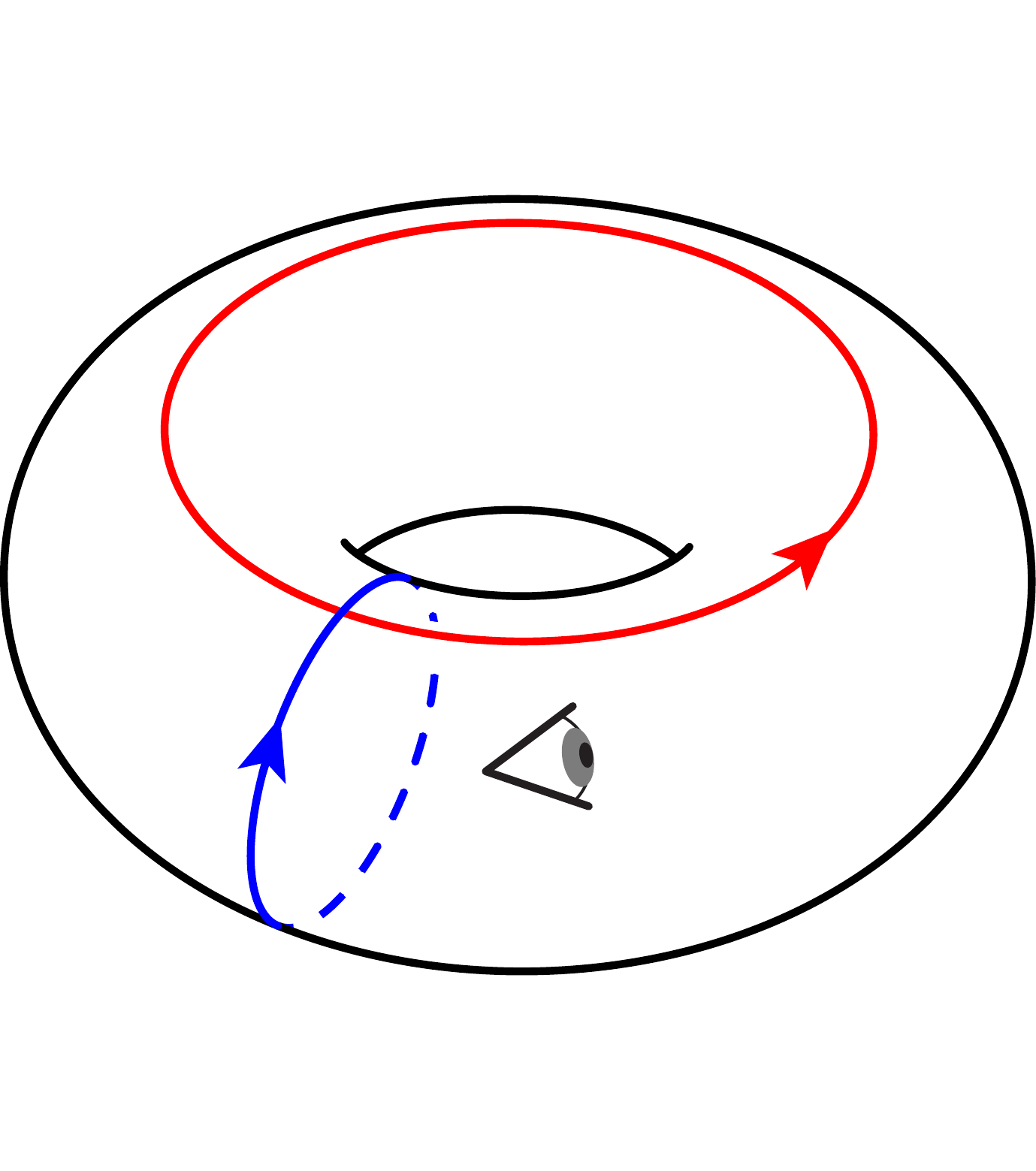}
}
\quad
\subfloat[\emph{The unrolled torus is a square.}]{
\label{Fig:TorusSquare}
\includegraphics[height=1.67in]{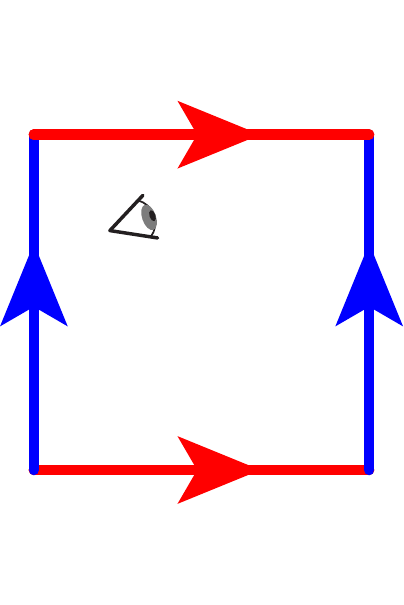}
}
\quad
\subfloat[\emph{The universal cover of the torus is the plane.}]{
\label{Fig:TorusUnivCover}
\includegraphics[height=1.67in]{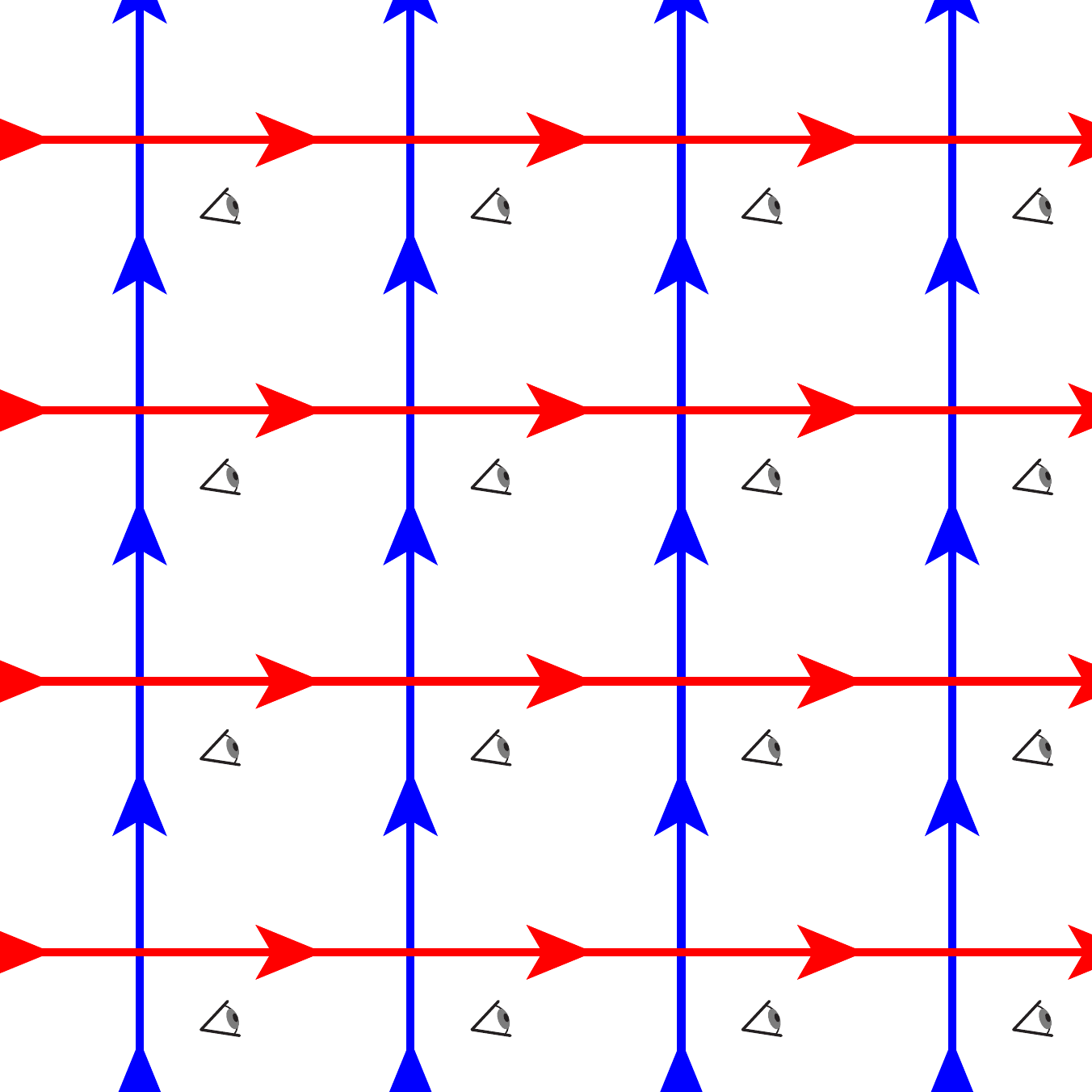}
}
\quad
\subfloat[\emph{The visual circle in the universal cover.}]{
\label{Fig:VisualCircle}
\includegraphics[height=1.5in]{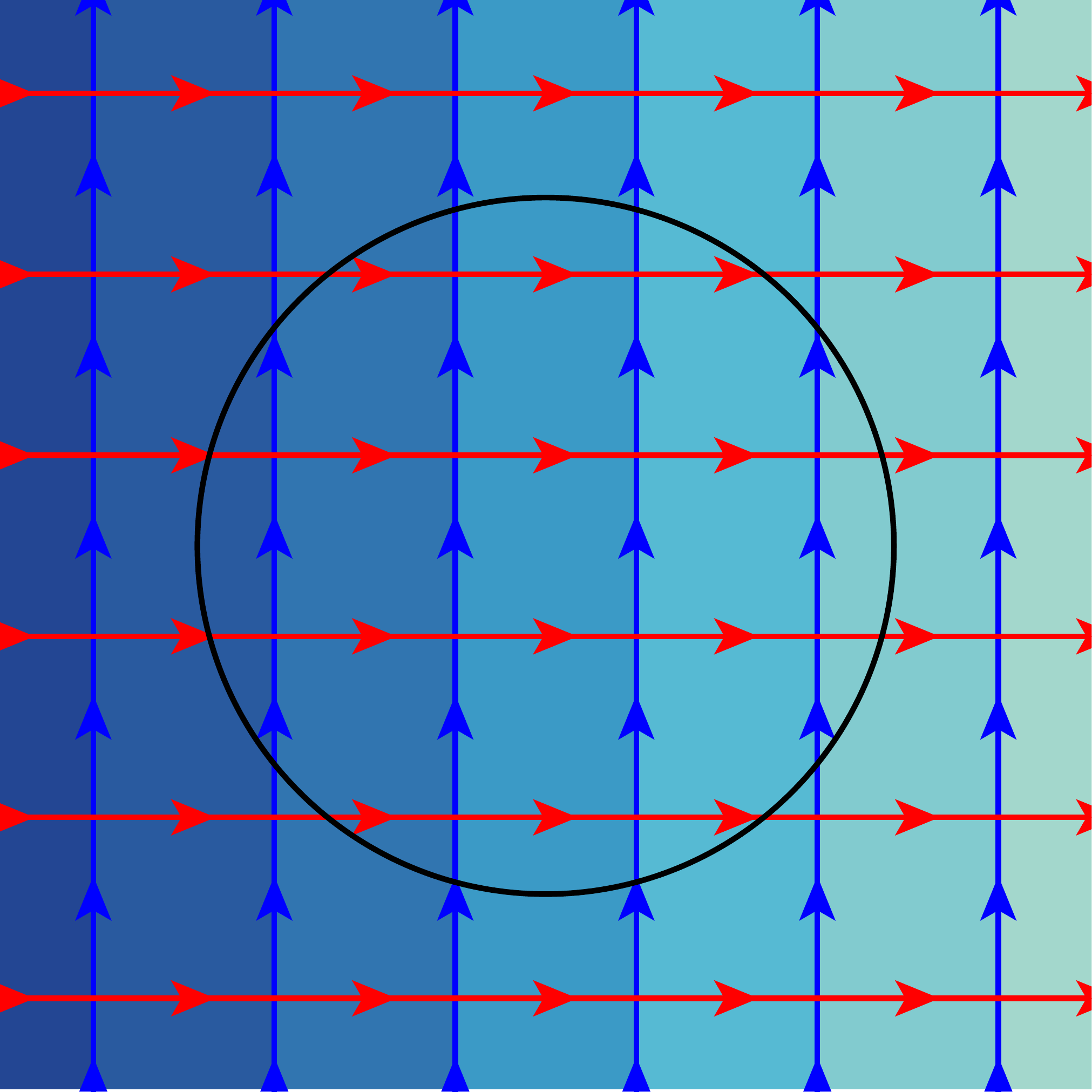}
}
\caption{}
\label{Fig:TorusAndCover}
\end{figure}

\section*{Warm-up examples: the torus and three-torus}

Before dealing with the three-dimensional hyperbolic case, we discuss the relevant ideas in simpler settings.  We refer to~\cite{Weeks19} for a richly illustrated introduction to surfaces and three-manifolds.  Consider the two-dimensional torus, as shown in \reffig{Torus}.  We have decorated the torus with two oriented circles: the blue circle goes through the hole of the torus, while the red circle goes around the hole.  Cutting along the blue circle turns the torus into a cylinder.  Cutting along the resulting red arc turns the cylinder into a square.  See \reffig{TorusSquare}.  Tiling this square out onto the plane produces the \emph{universal cover} of the torus, as shown in \reffig{TorusUnivCover}; see also~\cite[Section~1.3]{Hatcher02}.

In the universal cover, the red circle becomes an infinite collection of parallel horizontal lines, while the blue circle becomes an infinite collection of parallel vertical lines.  Each of these lines is called an \emph{elevation} of the corresponding circle; see~\cite[page~41]{Wise12} for a formal definition.

We now imagine ourselves living in the torus.  We suppose that light rays travel along the surface of the torus.  Thus, looking along the red curve, we see the back of our head.  This is because light rays bounce off our head, travel all the way around the torus, and then enter our eye.  The same is true when looking in the blue direction, or more generally, any direction with a rational slope.  So, we see a grid of copies of ourselves; this is the same as what we would see if we were living in the universal cover, tiled by squares.   However, in the universal cover the topology has been ``unwrapped'' and so the geometry is simpler; light here travels along straight lines.  Again, see \reffig{TorusUnivCover}.

Now suppose that the blue circle is a magic ``one-way'' filter, with the property that if we look through it in the direction of the red arrow things look lighter.  If we look through it in the opposite direction things look darker.  \reffig{VisualCircle} gives a ``birds-eye'' view of the universal cover of the torus, shaded according to these rules.  The further we look through the blue circle in the direction of the red arrow, the lighter it gets.  In \reffig{VisualCircle} we have also added a circle, in black.  If we look out to only that distance, then we see five different shades of light and dark.  We call this circle the \emph{visual circle} of a given radius -- in this example the radius is about twice the side length of the squares.  As we increase the radius of the visual circle, we see ever lighter and ever darker shades.  In the limit, a full half of our visual circle will be white, and the other half will be black.

\begin{figure}[htb]
\centering
\subfloat[\emph{The unrolled three-torus is a cube.}]{
\label{Fig:ThreeTorusCube}
\includegraphics[height=2in]{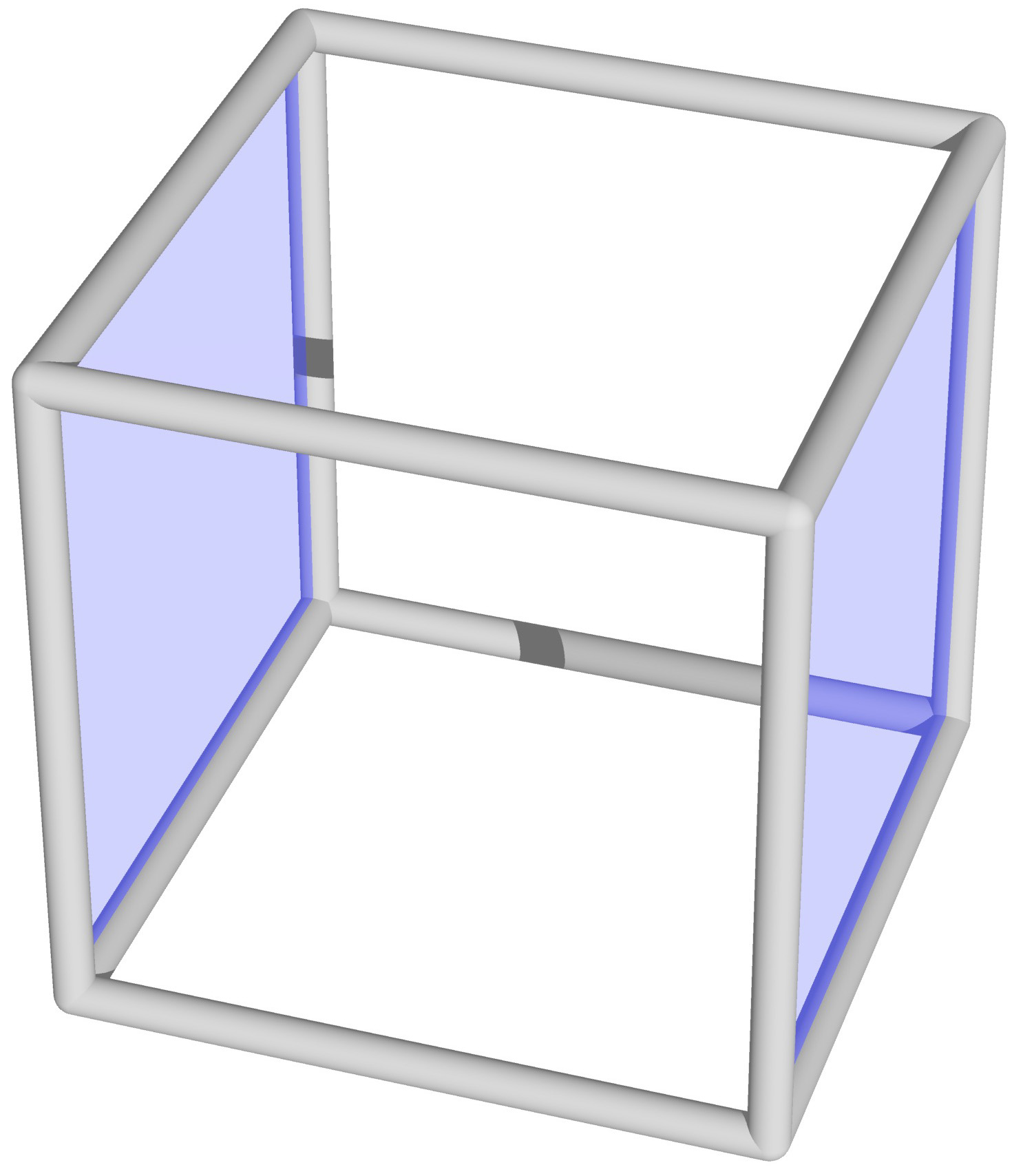}
}
\quad
\subfloat[\emph{The universal cover of the three-torus is three-dimensional space.}]{
\label{Fig:ThreeTorusUnivCover}
\includegraphics[height=2in]{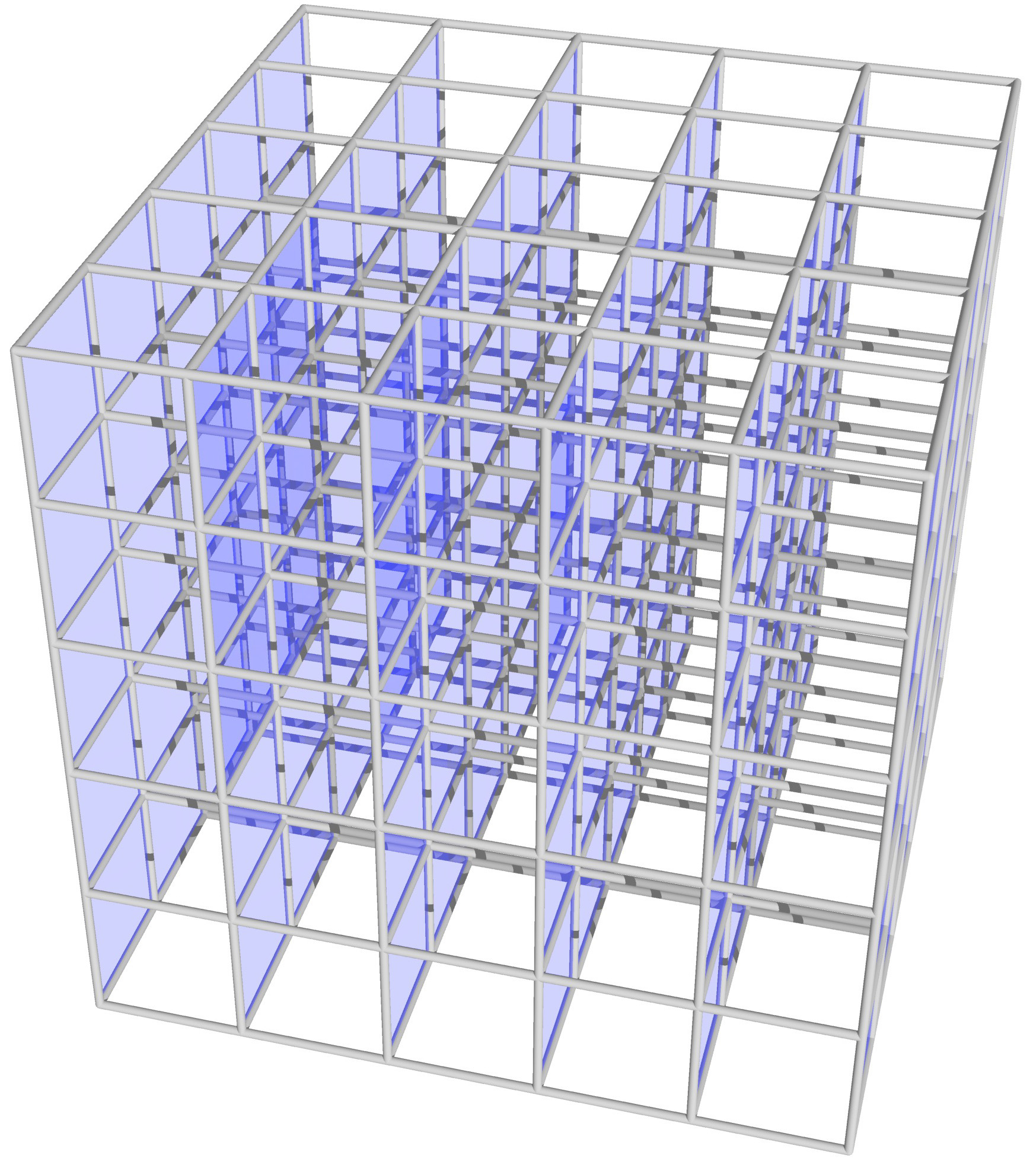}
}
\quad
\subfloat[\emph{The visual sphere in the universal cover.}]{
\label{Fig:VisualSphere}
\includegraphics[height=2in]{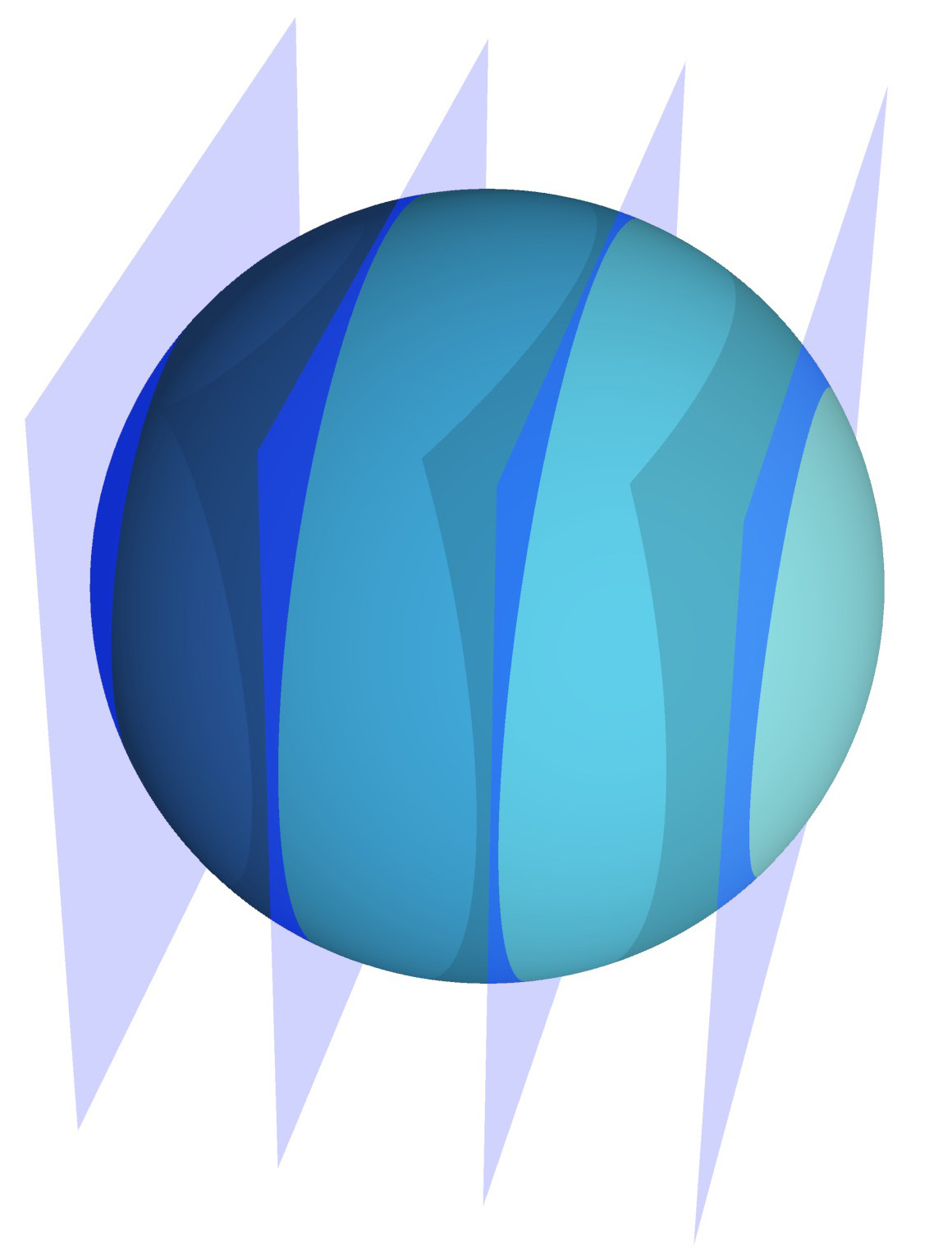}
}
\caption{}
\label{Fig:ThreeTorusAndCover}
\end{figure}

We now consider the three-dimensional torus~\cite[page~19]{Weeks19}.  Recall that we can build the torus shown in \reffig{Torus} by starting with the square in \reffig{TorusSquare} then gluing opposite sides to each other.  To construct the \emph{three-torus} we start with a cube, then glue opposite sides to each other.  See \reffig{ThreeTorusCube}. 

We colour, in blue, a pair of parallel walls of the cube.  These glue up to give a two-torus inside of the three-torus.  Each elevation of the blue two-torus is a plane in the universal cover of the three-torus.  See \reffig{ThreeTorusUnivCover}.  The visual circle is replaced by the \emph{visual sphere}; again one side of it will be lighter, and the other will be darker.  See \reffig{VisualSphere}.  As the radius of the sphere increases to infinity, a hemisphere becomes pure white while the opposite hemisphere becomes pure black.

\section*{Constructing cohomology fractals}

\begin{figure}[htb]
\centering
\subfloat[\emph{Tube around the figure-eight knot.}]{
\label{Fig:Fig8Knot}
\includegraphics[height=2in]{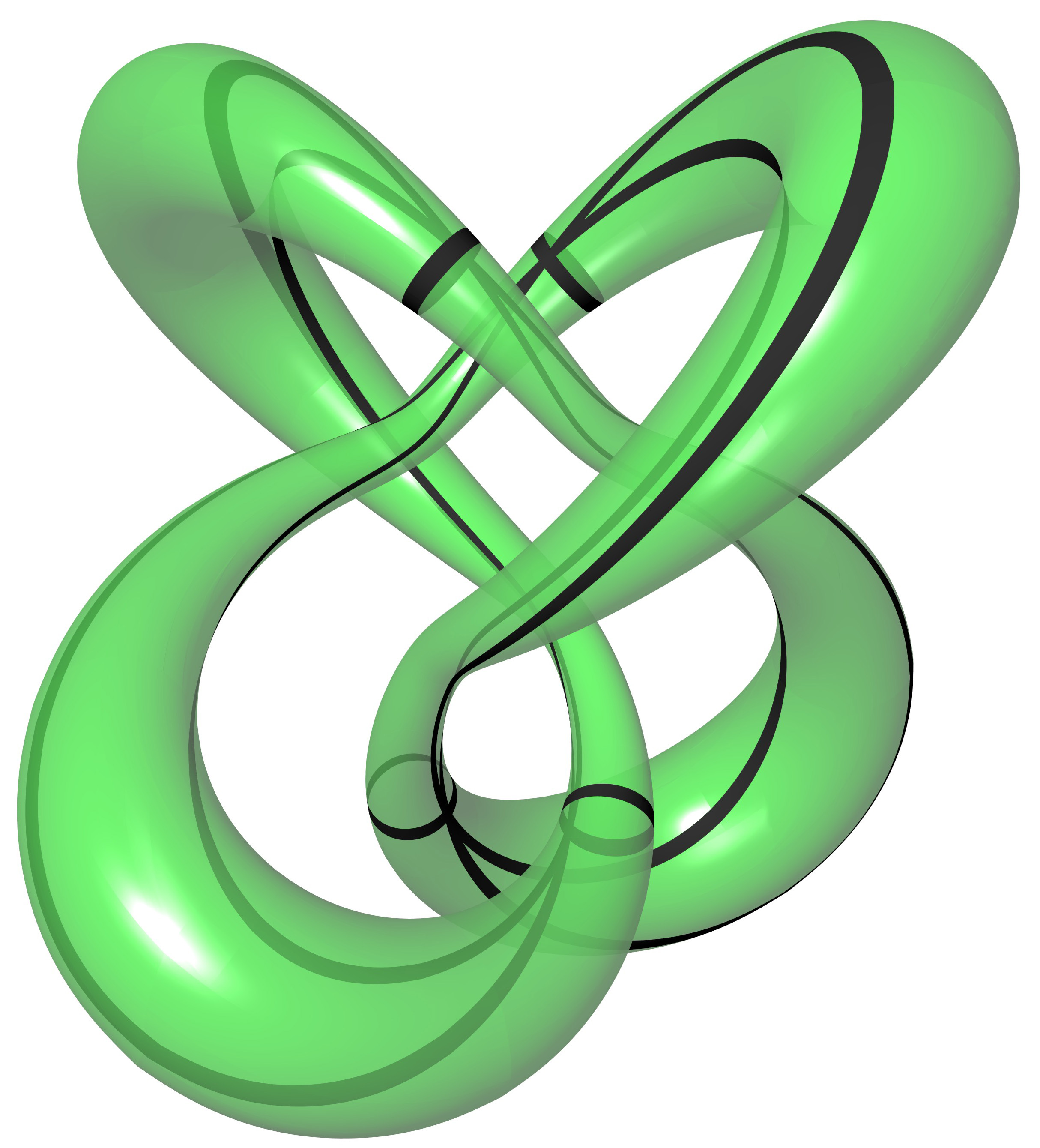}
}
\qquad
\subfloat[\emph{Triangulation of the figure-eight knot complement.}]{
\label{Fig:Fig8Triangulation}
\includegraphics[height=1.9in]{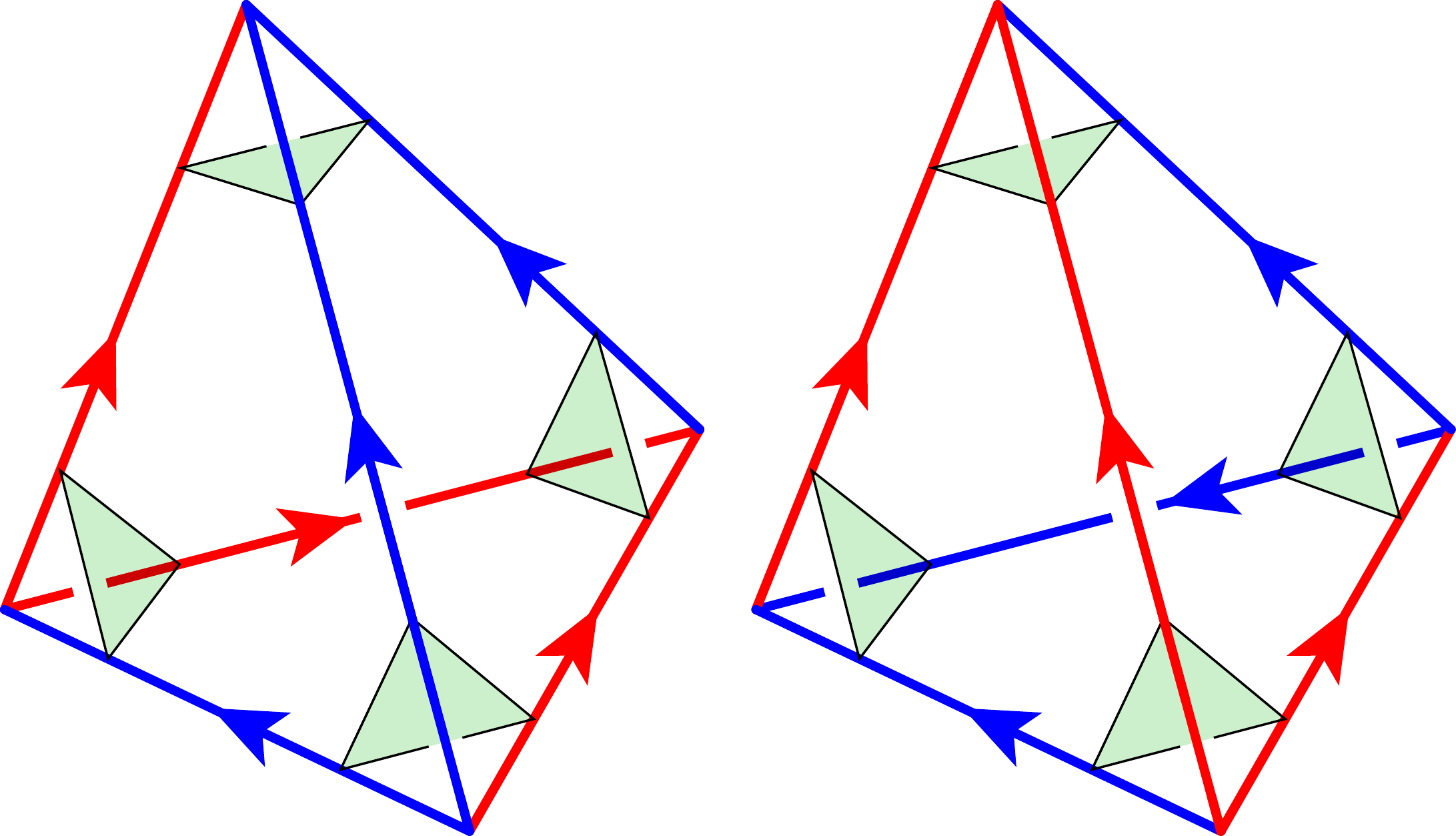}
}

\subfloat[\emph{A tiling of $\HH^3$ by ideal hyperbolic tetrahedra.}]{
\label{Fig:H3Tiling}
\includegraphics[width=0.44\textwidth]{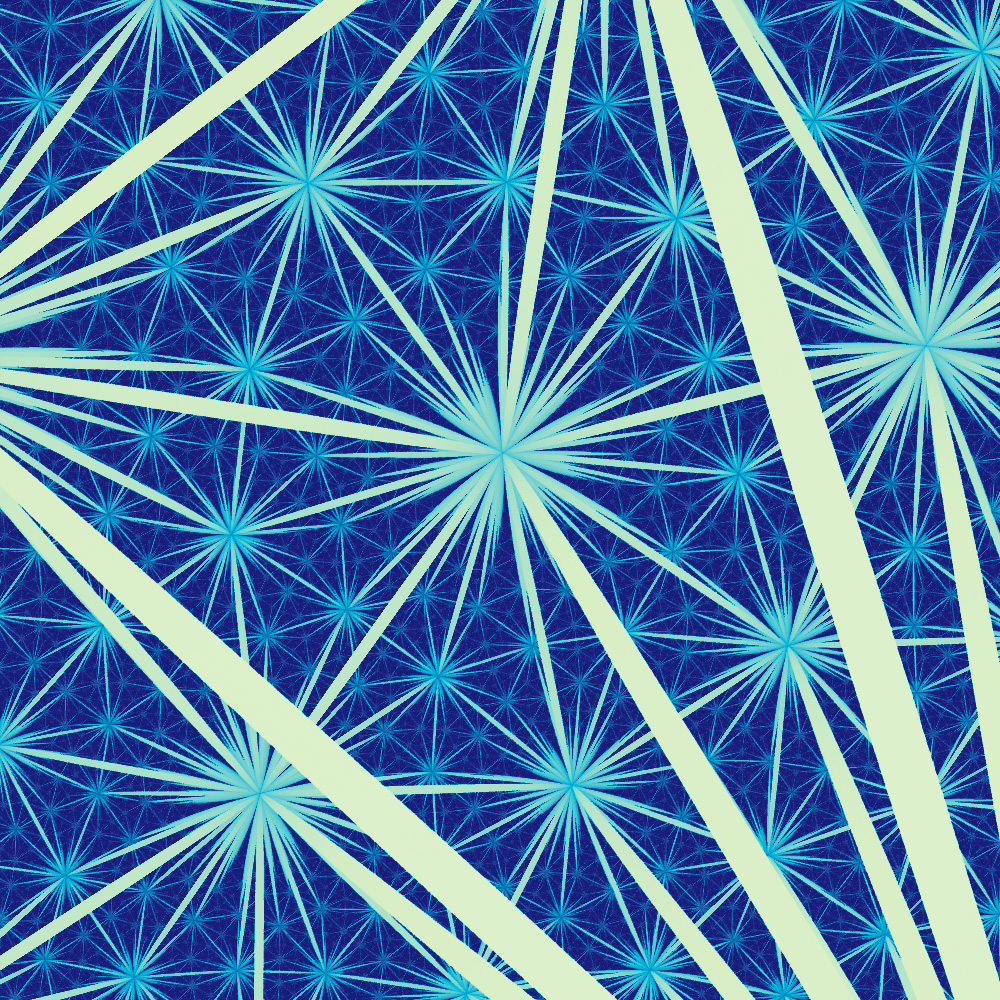}
}
\quad
\subfloat[\emph{An elevation of the surface inside of $\HH^3$.}]{
\label{Fig:Elevation}
\includegraphics[width=0.44\textwidth]{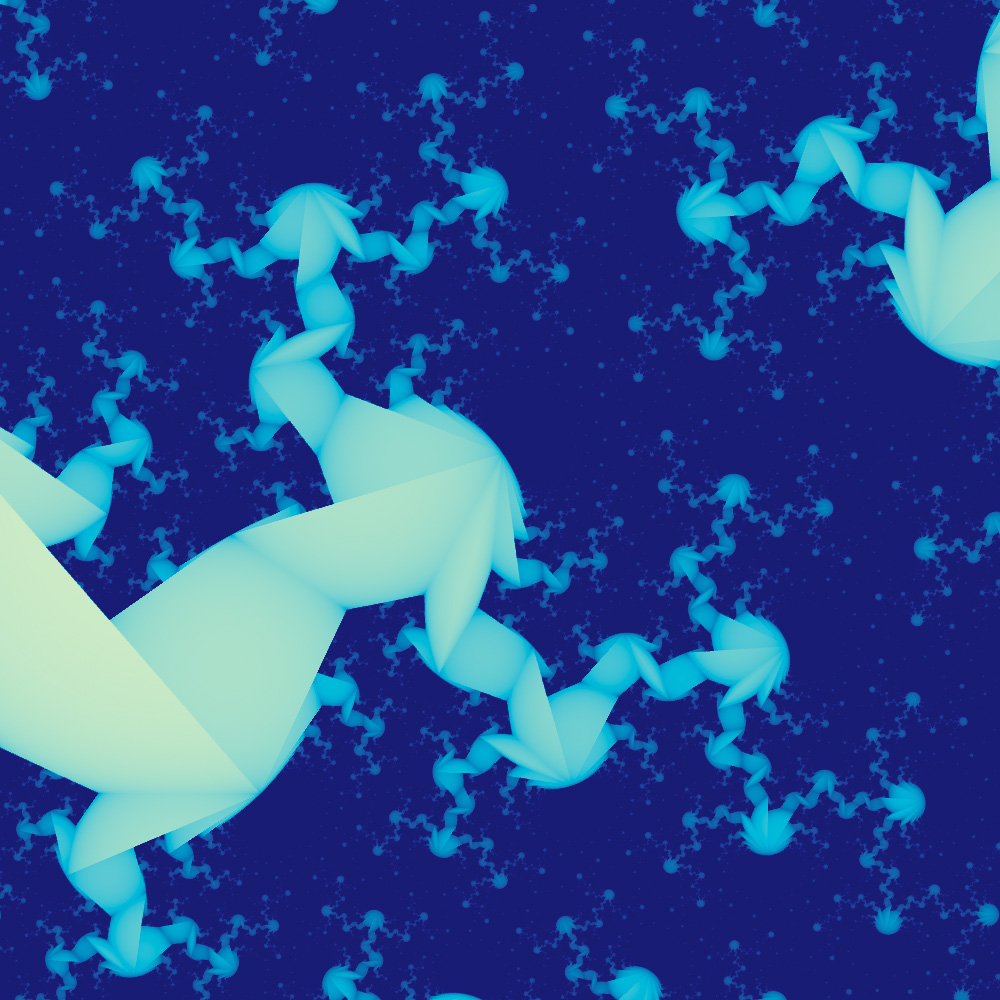}
}
\caption{\emph{The figure-eight knot complement.}}
\label{Fig:Fig8KnotExample}
\end{figure}

We produce our images in much the same way as the previous example, but with one key difference: we replace the euclidean three-torus with a three-dimensional \emph{hyperbolic} manifold $M$ with torus boundary components.  We refer to~\cite[Chapter~2]{PurcellKnotTheory} for an introduction to hyperbolic geometry.

The universal cover of the three-torus is a euclidean space, tiled by cubes; again see~\reffig{ThreeTorusUnivCover}.  The universal cover of $M$ is three-dimensional hyperbolic space $\HH^3$, tiled by \emph{ideal tetrahedra}~\cite[Definition~2.17]{PurcellKnotTheory}.  Here an ideal hyperbolic tetrahedron sits inside of $\HH^3$ with each face a subset of a geodesic plane, each edge a geodesic, and each vertex a point on the ``sphere at infinity'' $\bdy \HH^3$.  We draw this tiling in \reffig{H3Tiling} as viewed from inside of $\HH^3$.  Here we again assume that light rays travel along geodesics. See~\cite{visualizing_hyperbolic_honeycombs} for more on tilings of hyperbolic space.

The manifold $M$, shown in \reffig{Fig8KnotExample}, is the \emph{complement of the figure-eight knot} in the \emph{three-sphere}~\cite[Section~9.2]{Adams04}; in the \texttt{SnapPy} census~\cite{snappy} it is called \texttt{m004}.  We delete the figure-eight knot (shown in \reffig{Fig8Knot}) from the three-sphere, and what remains is the manifold $M$.  William Thurston~\cite[Chapter~1, page~4]{ThurstonNotes} showed that $M$ can be subdivided into two ideal tetrahedra.  Jessica Purcell gives an excellent exposition of Thurston's construction in~\cite[Section~2.3]{PurcellKnotTheory}.  The third author has made a sculpture illustrating how the two tetrahedra twist to fit around the figure-eight knot~\cite[page~141]{3dprintmath}.  

This triangulation is illustrated in \reffig{Fig8Triangulation}.  The four faces of the first tetrahedron are glued to the four faces of the second tetrahedron in such a way that the arrows and colours on the edges match.  These also induce gluings of green triangles near to the vertices of the tetrahedra together.  These triangles join to form a torus around the ``missing'' knot; see \reffig{Fig8Knot}.  

For our one-way filter, we take take two of the faces of the tetrahedra.  Here we choose the ones with exactly one blue and two red edges; see \reffig{Fig8Triangulation}.  This is in fact a \emph{Seifert surface}~\cite[page~99]{Adams04} for the figure-eight knot.  One elevation of the surface to the universal cover is shown in \reffig{Elevation}.  Unlike the elevation in the three-torus, now it is an intricately pleated surface. 

\begin{figure}[htb!]
\centering
\subfloat[$R=e^0$]{
\label{Fig:exp0}
\includegraphics[width=0.32\textwidth]{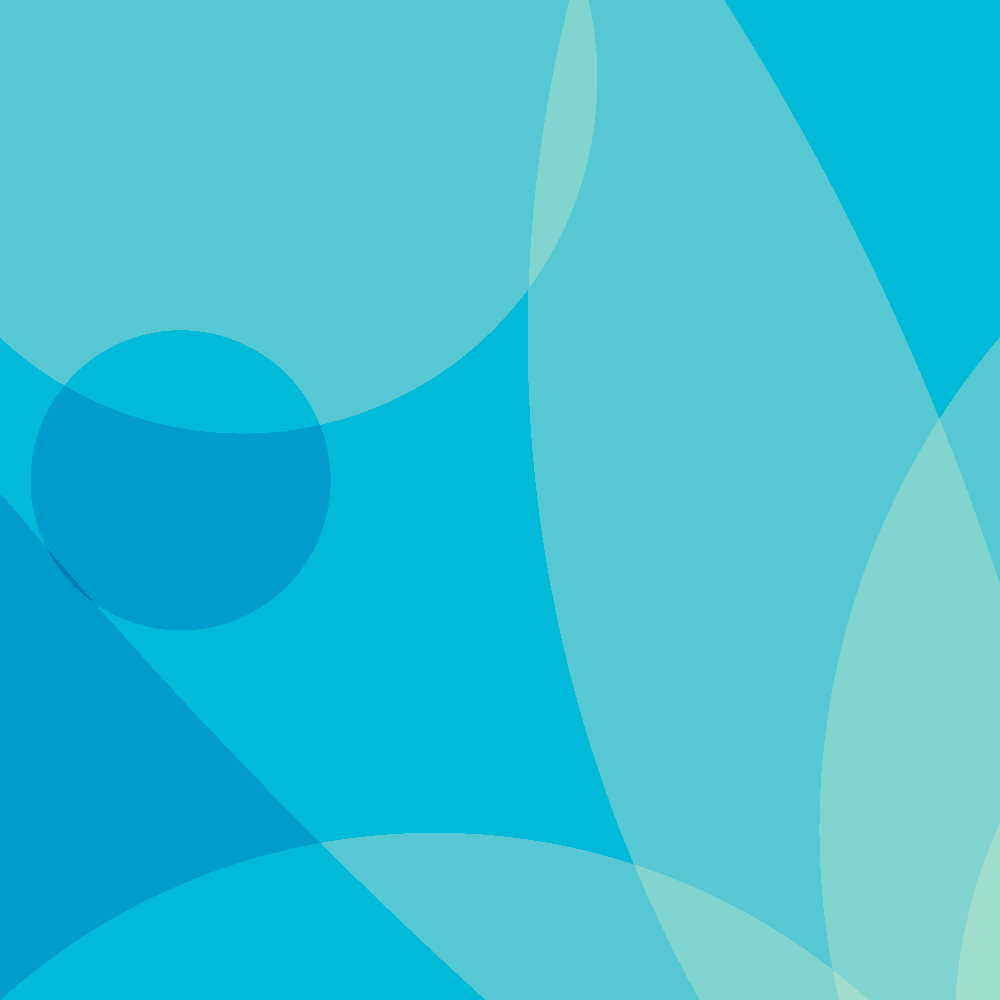}
}
\subfloat[$R=e^{0.5}$]{
\label{Fig:exp0.5}
\includegraphics[width=0.32\textwidth]{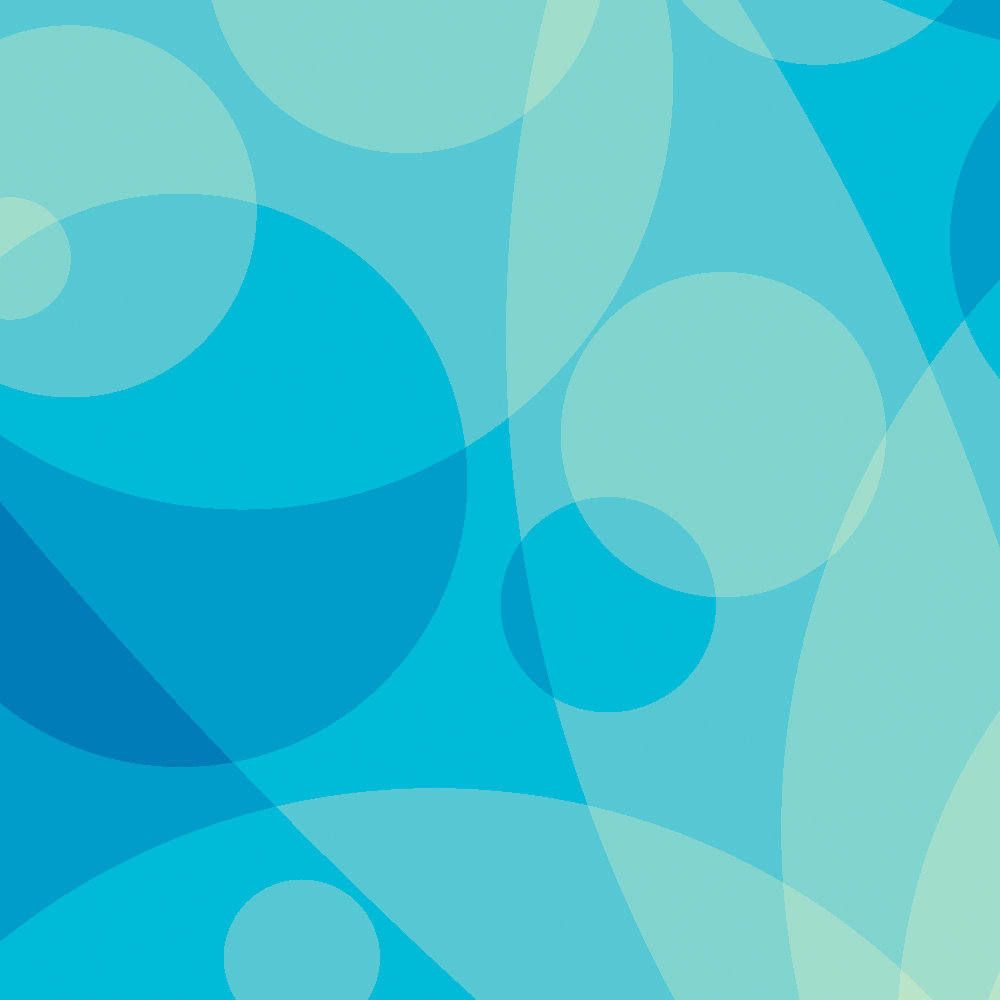}
}
\subfloat[$R=e^1$]{
\label{Fig:exp1}
\includegraphics[width=0.32\textwidth]{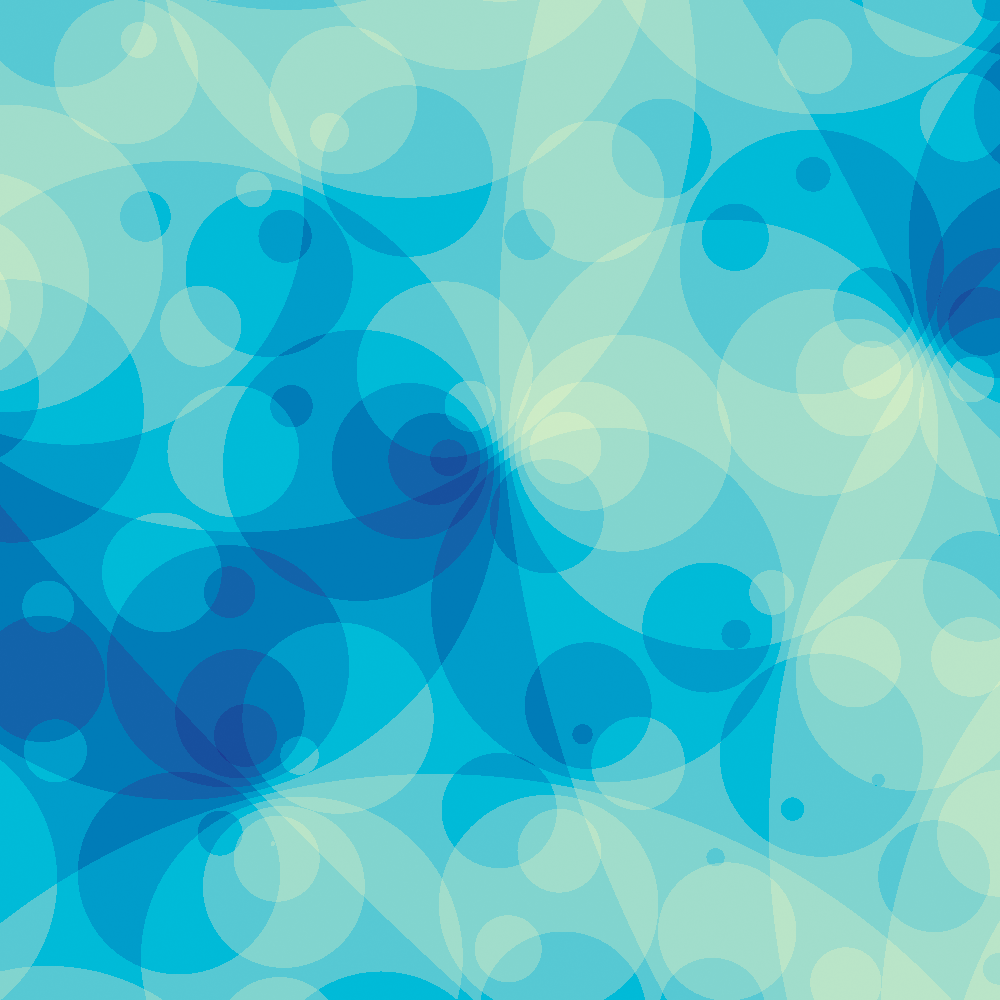}
}

\subfloat[$R=e^{1.5}$]{
\label{Fig:exp1.5}
\includegraphics[width=0.32\textwidth]{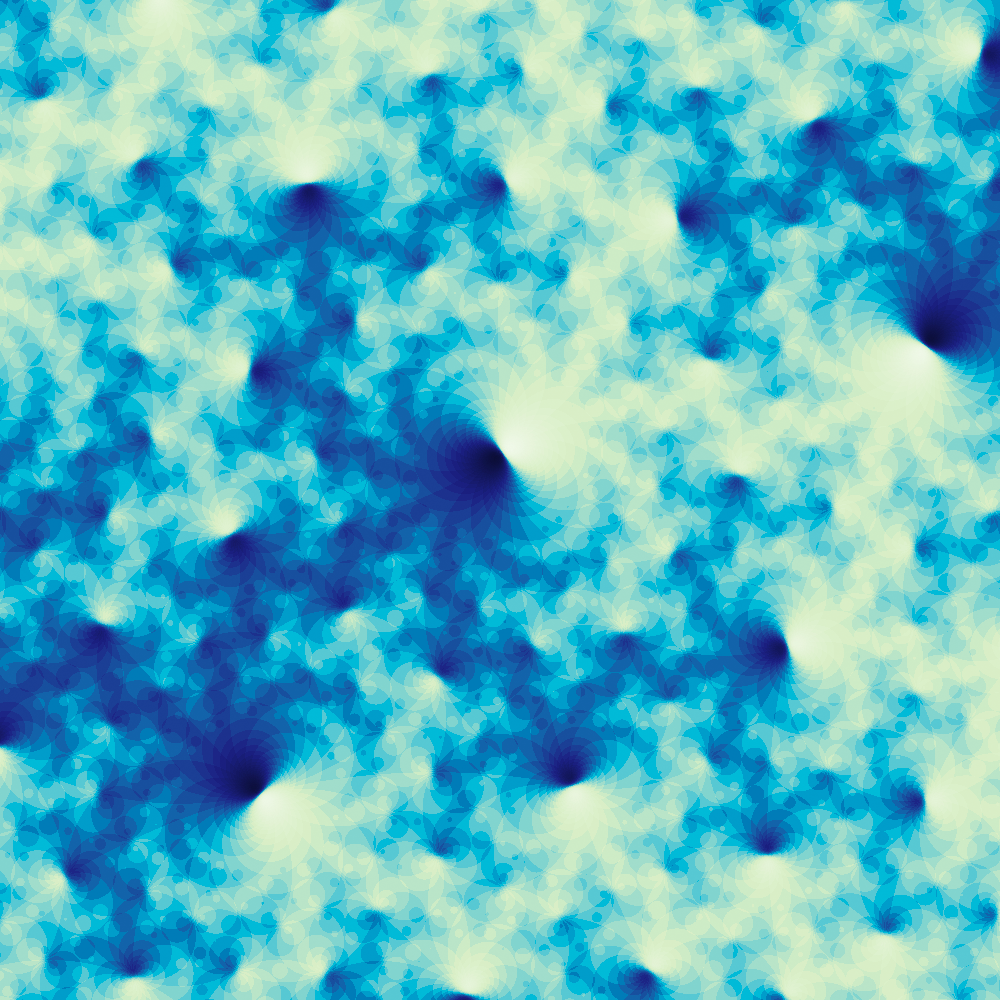}
}
\subfloat[$R=e^{2}$]{
\label{Fig:exp2}
\includegraphics[width=0.32\textwidth]{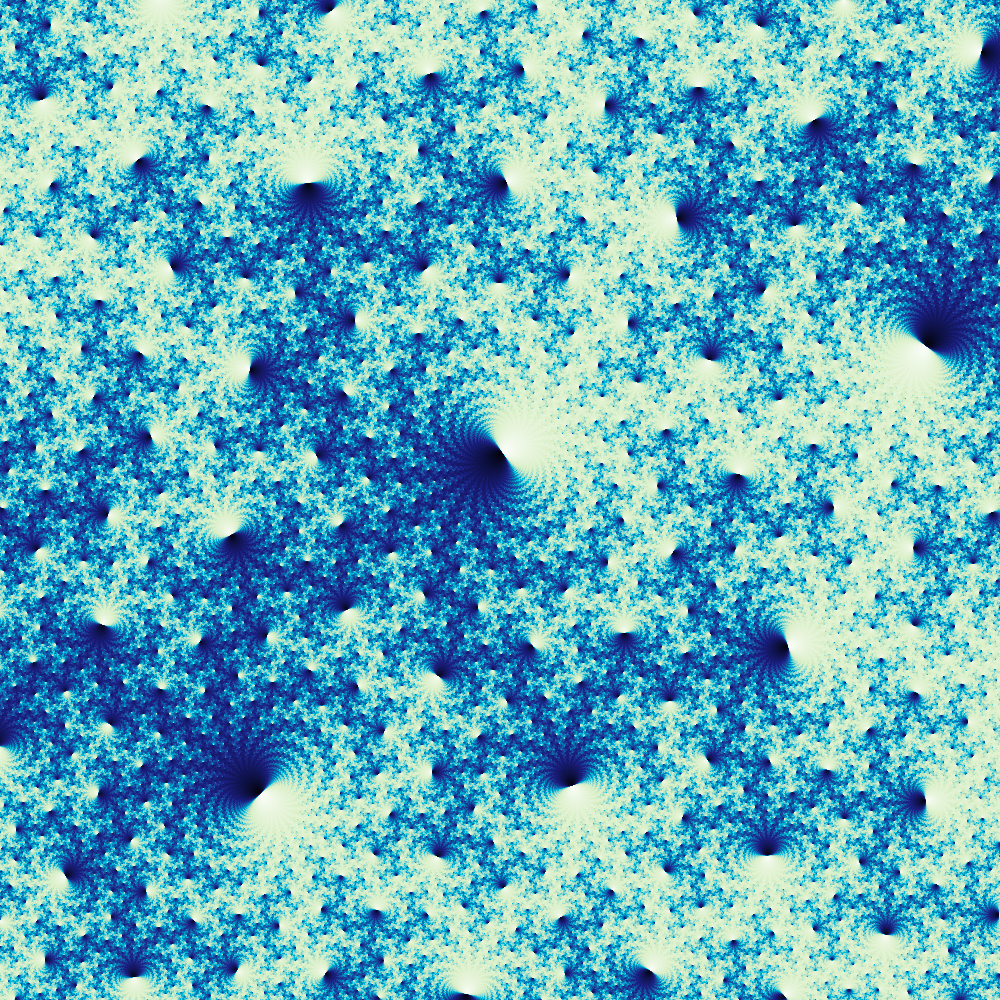}
}
\subfloat[$R=e^{2.5}$]{
\label{Fig:exp2.5}
\includegraphics[width=0.32\textwidth]{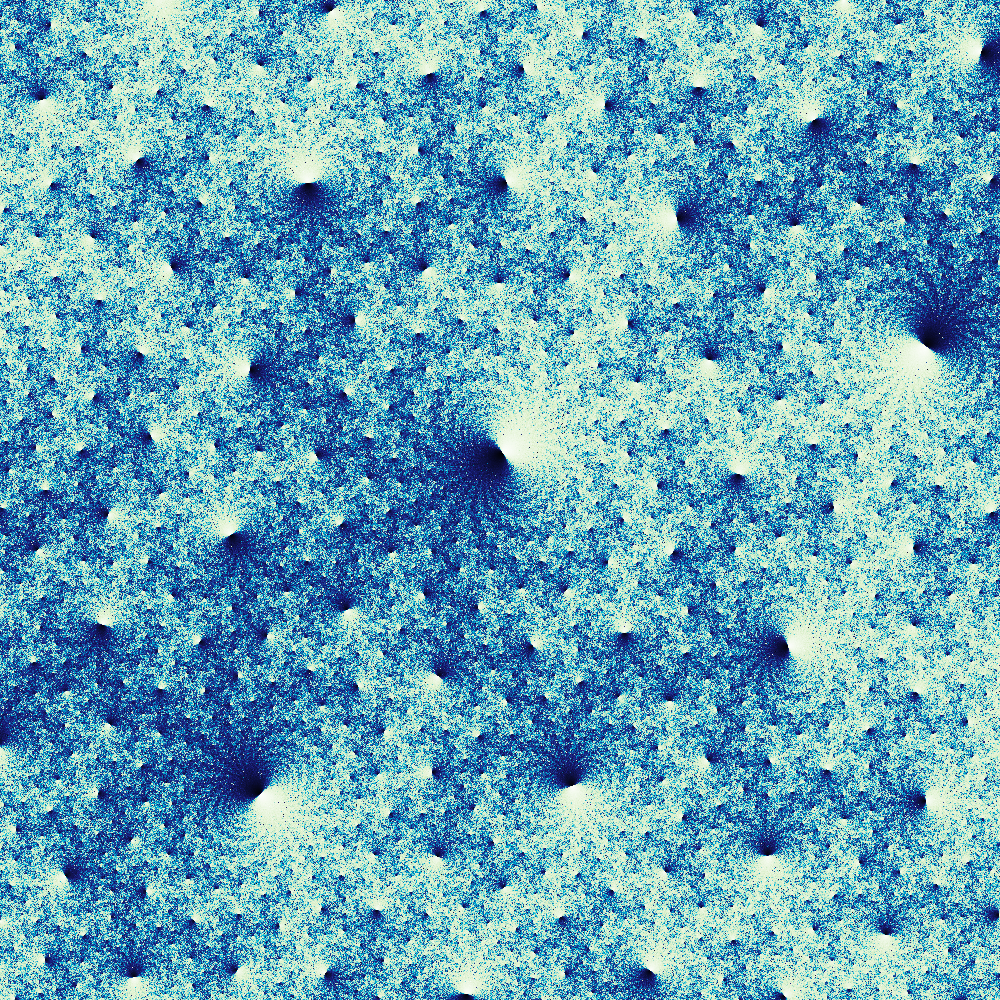}
}

\subfloat[$R=e^{3}$]{
\label{Fig:exp3}
\includegraphics[width=0.32\textwidth]{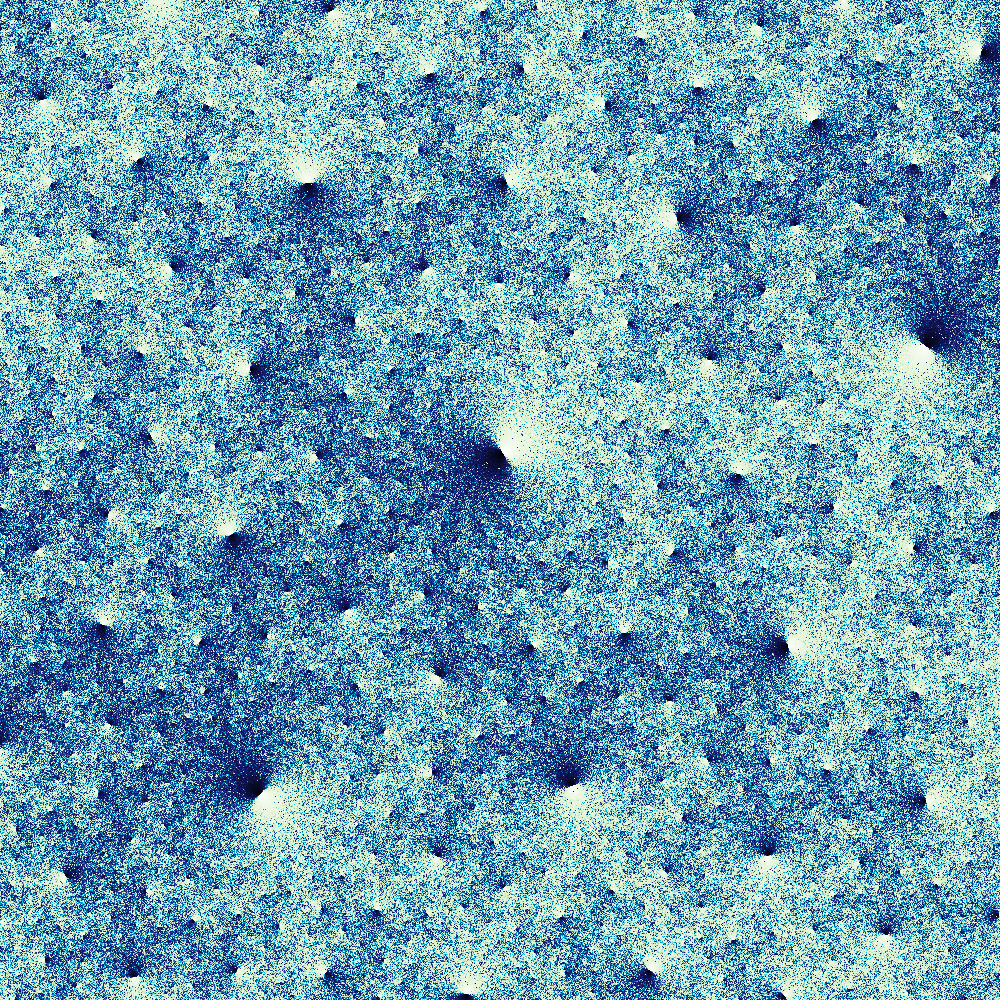}
}
\subfloat[$R=e^{4}$]{
\label{Fig:exp4}
\includegraphics[width=0.32\textwidth]{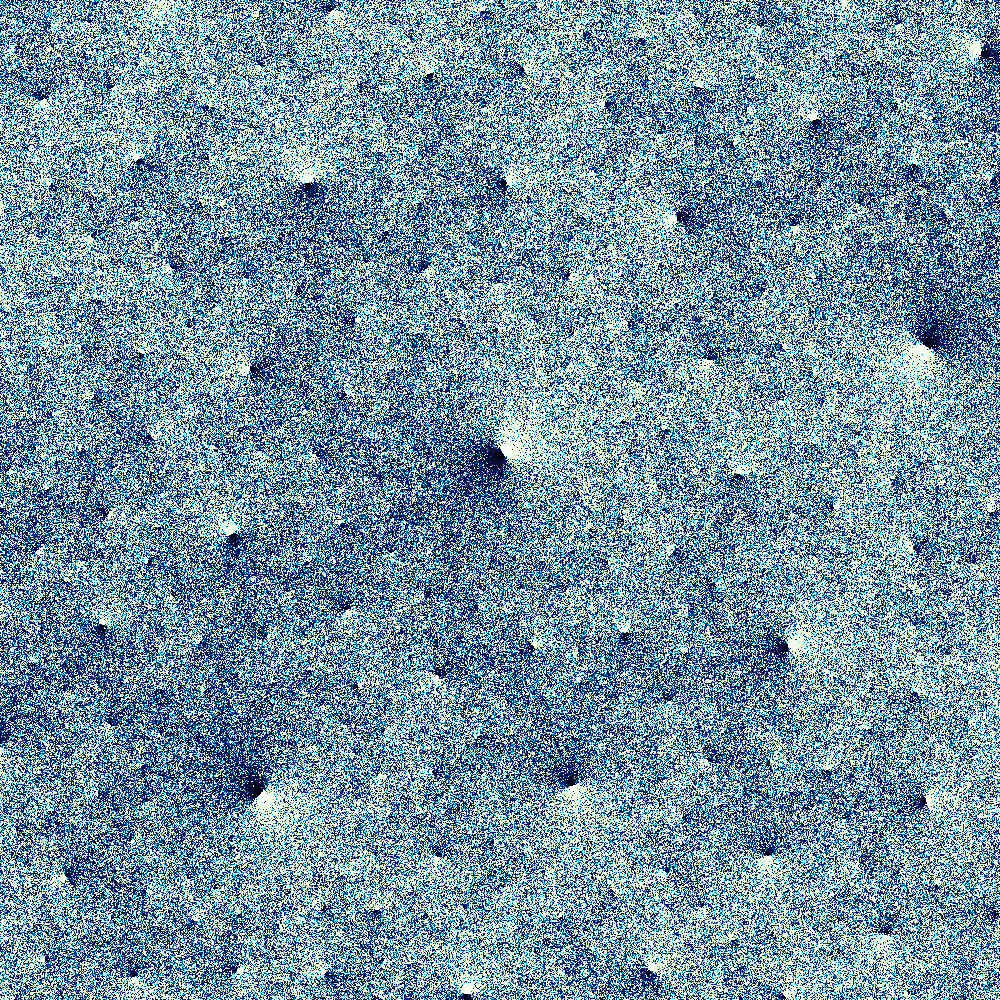}
}
\subfloat[$R=e^{5}$]{
\label{Fig:exp5}
\includegraphics[width=0.32\textwidth]{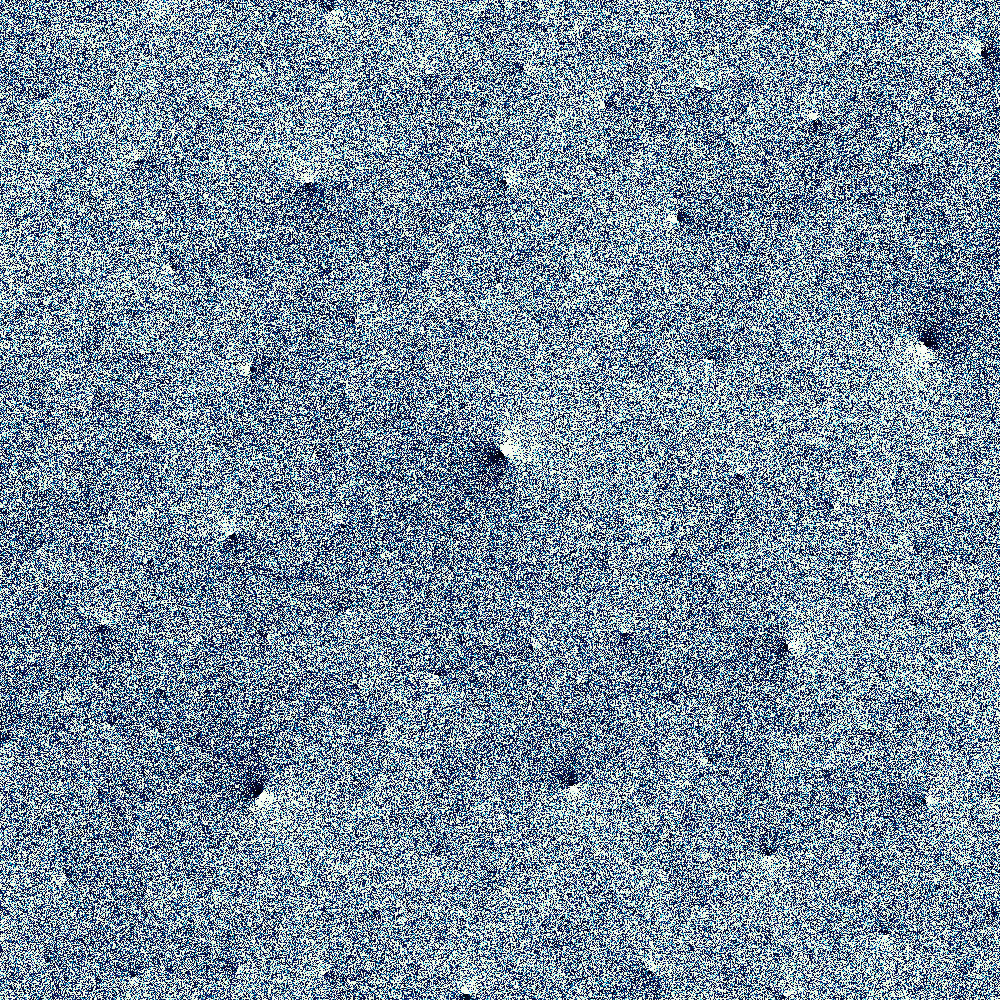}
}
\caption{\emph{Cohomology fractals for} \texttt{m004}, \emph{with various values of $R$.}}
\label{Fig:VisSphereRadii}
\end{figure}

Again, we imagine ourselves in the manifold $M$, looking along geodesics, with all filter surfaces in place.  \reffig{VisSphereRadii} shows the resulting view with various choices of radius $R$ for the visual sphere, but always with a fixed number of pixels.
These then are the \emph{cohomology fractals}.  (We will explain the relevance of cohomology in the section on implementation.)  Certainly, \reffig{exp2} has a fractal appearance.  When $R$ is small, pixels are tightly correlated with their neighbours.  When $R$ is large, they are not.  Images in the middle of the sequence, like \reffig{exp2}, have the property that each pixel (at that radius) meets approximately one tetrahedron.  If we increase the number of pixels by a factor of $k^2$, and simultaneously increase the radius of the visual sphere by $\log(k)$, then the sharpness of the image remains constant.

\section*{Recurring motifs: lighthouses, blind spots, circles, spirals}

The cohomology fractals of \reffig{PixelImages} have some common features; we explain some of these in terms of the hyperbolic geometry and/or the topology of the manifold $M$ and the surface $F$ it contains. 

In all of the images we see \emph{lighthouses}: that is, points equipped with a tangent line where the image is very bright on one side of the line and very dark on the other.  
Each lighthouse lies at a vertex of some ideal hyperbolic tetrahedron in $\cover{M}$, the universal cover of $M$.  A neighbourhood of an ideal point is called a \emph{cusp neighbourhood}~\cite[Definition~4.9]{PurcellKnotTheory}.  In $M$ the cusp neighbourhood is a copy of a torus crossed with a ray; in $\cover{M}$ the cusp neighbourhood is a copy of an open ball.  If one elevation of the surface $F$ enters a given cusp neighbourhood, then so do infinitely more.  Furthermore all elevations entering the cusp are essentially parallel.  When a ray enters the cusp neighbourhood, it may cross many elevations very quickly.  This gives the large weights, and thus the very bright and very dark pixels on either side of the lighthouse. 

Most of the three-manifolds used in \reffig{PixelImages} have only one cusp.  We nonetheless see infinitely many lighthouses in each image.  Again, this is because sitting inside the manifold we see its universal cover.  Seeing features from different distances and angles contributes to the fractal appearance of the image. 

The manifold \texttt{m129}, also known as the \emph{Whitehead link complement}~\cite[Example~6.16]{PurcellKnotTheory}, has two cusps.  It follows that \texttt{m129} has infinitely many distinct cohomology fractals.  
In \reffig{m129} we show its simplest cohomology fractal; here the surface is a once-holed torus $T$.  Let $N(c)$ and $N(d)$ be cusp neighbourhoods of the cusps $c$ and $d$ of the manifold.  Since $T$ has only one boundary component we may assume that it meets $N(c)$ and not $N(d)$.  So rays entering $N(c)$ can meet elevations of $T$; this leads to lighthouses at all of the lifts of $c$.  On the other hand, rays cannot meet $T$ while inside of $N(d)$; this leads to \emph{blind spots} at all lifts of $d$.  One of these is shown at the exact centre of \reffig{m129}.  Rays entering a cusp neighbourhood have commuting choices for how they can wind about the two directions of the torus.  Again, see \reffig{Torus}.  Thus there is (an inversion of) a copy of the universal cover of the torus (as in \reffig{TorusUnivCover}) packed into a grid-like pattern about each cusp of $\cover{M}$. 

\begin{figure}[htbp]
\centering
\subfloat[\emph{Super-apollonian packing.}]{
\label{Fig:Apollonian}
\includegraphics[width=0.32\textwidth]{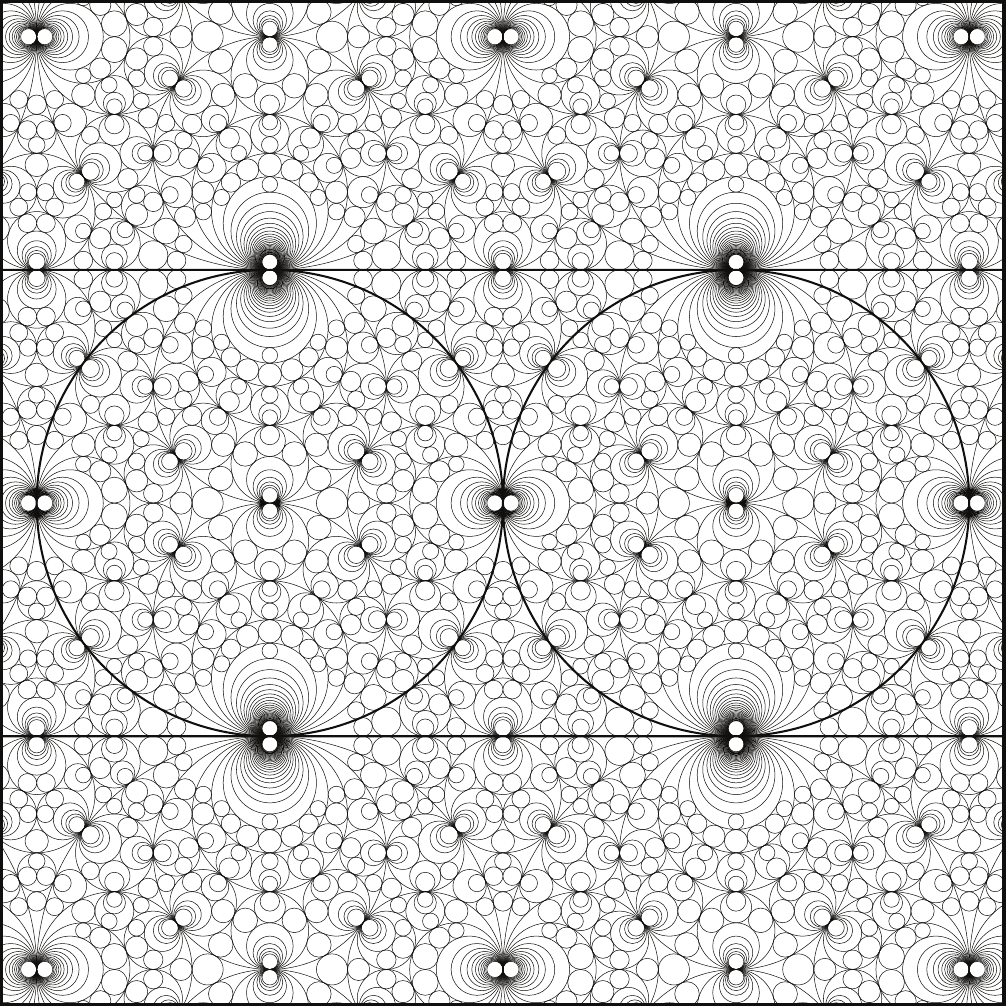}
}
\subfloat[\emph{Cohomology fractal for} \texttt{m129}.]{
\label{Fig:m129Overlay}
\includegraphics[width=0.32\textwidth]{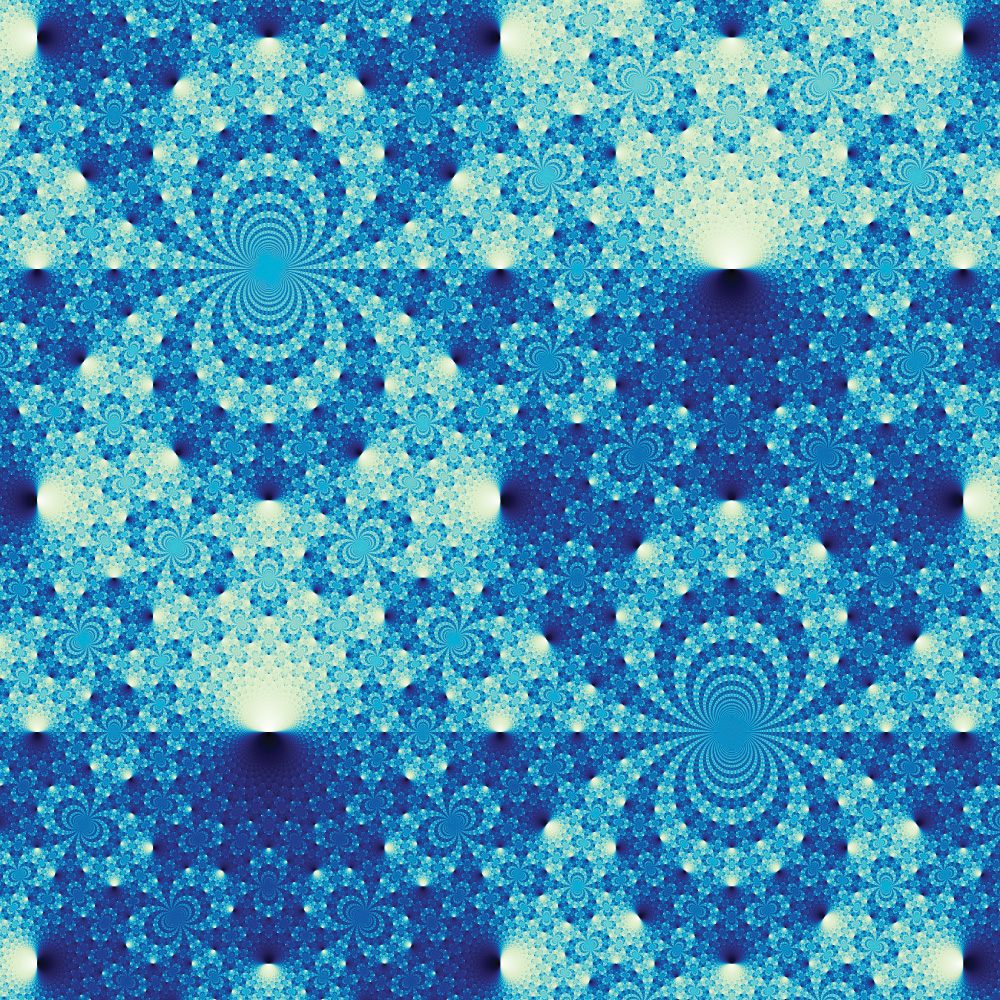}
}
\subfloat[\emph{The super-apollonian packing in red, overlaid on the cohomology fractal.}]{
\label{Fig:Both}
\includegraphics[width=0.32\textwidth]{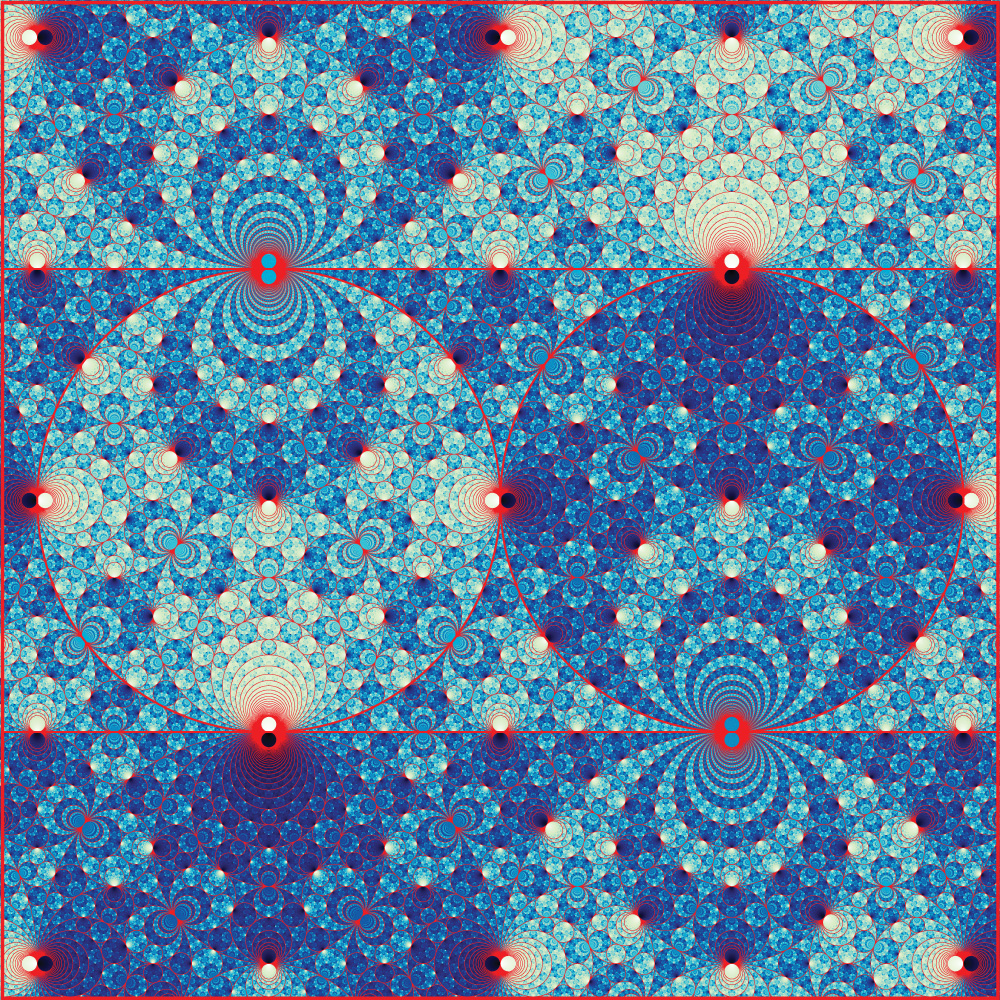}
}
\caption{\emph{Comparison between a super-apollonian circle packing}~\cite[Figure~4]{Graham06} \emph{and a cohomology fractal for the Whitehead link complement,} \texttt{m129}.}
\label{Fig:Graham2005_and_WH_link}
\end{figure}

The images for \texttt{m129} and \texttt{s776}, as shown in Figures~\ref{Fig:m129} and~\ref{Fig:s776}, consist of circles.  This is because the surfaces used are thrice-holed spheres (or become such after an annular compression).  Thrice-holed spheres, in hyperbolic geometry, always give totally geodesic planes~\cite[page~149]{ThurstonNotes}.  
These have circle boundary at infinity, explaining the appearance of the figures.  
In fact, \reffig{m129} contains a \emph{super-apollonian circle packing} in the sense of~\cite[Section~2]{Graham06}.  We give a ``proof-by-picture'' of this in \reffig{Graham2005_and_WH_link}.  \reffig{s776} does not contain a Descartes configuration of four circles, so it is not a super-apollonian circle packing. 

The spiralling features in \reffig{s000} come from \emph{loxodromic} isometries of hyperbolic space~\cite[Theorem~2.16]{PurcellKnotTheory}.  These correspond to the ``screw motions'' of euclidean geometry; a rotation in a plane followed by a translation perpendicular to that plane.  Given a loxodromic $\gamma$, we can change coordinates in the \emph{Poincar\'e ball model}~\cite[page~121]{Adams04} so that $\gamma$ fixes the north and south poles.  Then $\gamma$ moves a feature, say a lighthouse, along a \emph{loxodrome} (a curve of constant angle to the lines of longitude).  When the loxodrome is stereographically projected to the plane (as in our images) it gives an equiangular spiral. 

\section*{Implementation}
\label{Sec:Implement}

We obtain the hyperbolic three-manifolds from the software \texttt{SnapPy}~\cite{snappy}.  Each manifold $M$ is given by an ideal triangulation $\calT$, with a hyperbolic \emph{shape} for each ideal tetrahedron.  These glue together, via isometries, to give a geometric structure on $M$.  

Our one-way filter surfaces come from \emph{cohomology classes} for $M$~\cite[Chapter~3]{Hatcher02}.  Informally, the first cohomology group $H^1(M;\ZZ)$ consists of certain functions $\rho$ that take in oriented loops and return integers.  In practice, to represent $\rho$ we place weights on the (transversely oriented) faces of the triangulation; a loop picks up the weight, with sign, when it crosses a face in the positive or negative sense.  There are several ways to find weights that represent a cohomology class $\rho$.  A particularly nice (and finite) selection of classes is given by the vertices of the \emph{Thurston norm ball}~\cite[Theorem~1]{Thurston86}.  To compute these we relied on the software \texttt{tnorm}~\cite{tnorm20} and \texttt{regina}~\cite{regina}, as well as our own code. 


The key operation of our algorithm is hyperbolic \emph{ray-casting}.  We cast geodesic rays~\cite[Der Zeichner der Laute]{Durer25}
through the tetrahedra of the triangulation, summing the weights the ray picks up as it passes through faces.  For each pixel, we take the resulting integer, and apply the function $x \mapsto 1/2(1 + x/(|x|+1))$.  (This is graphically similar to $\arctan(x)$ but is computationally cheaper.)  We map the resulting value in $[0,1]$ to the colour of the pixel.  Rather than mapping directly to brightness, we go from black to white via a few other colours.  This allows the viewer to see more detail.

\begin{wrapfigure}[18]{r}{0.55\textwidth}
\vspace{15pt}
\centering
\includegraphics[width=0.53\textwidth]{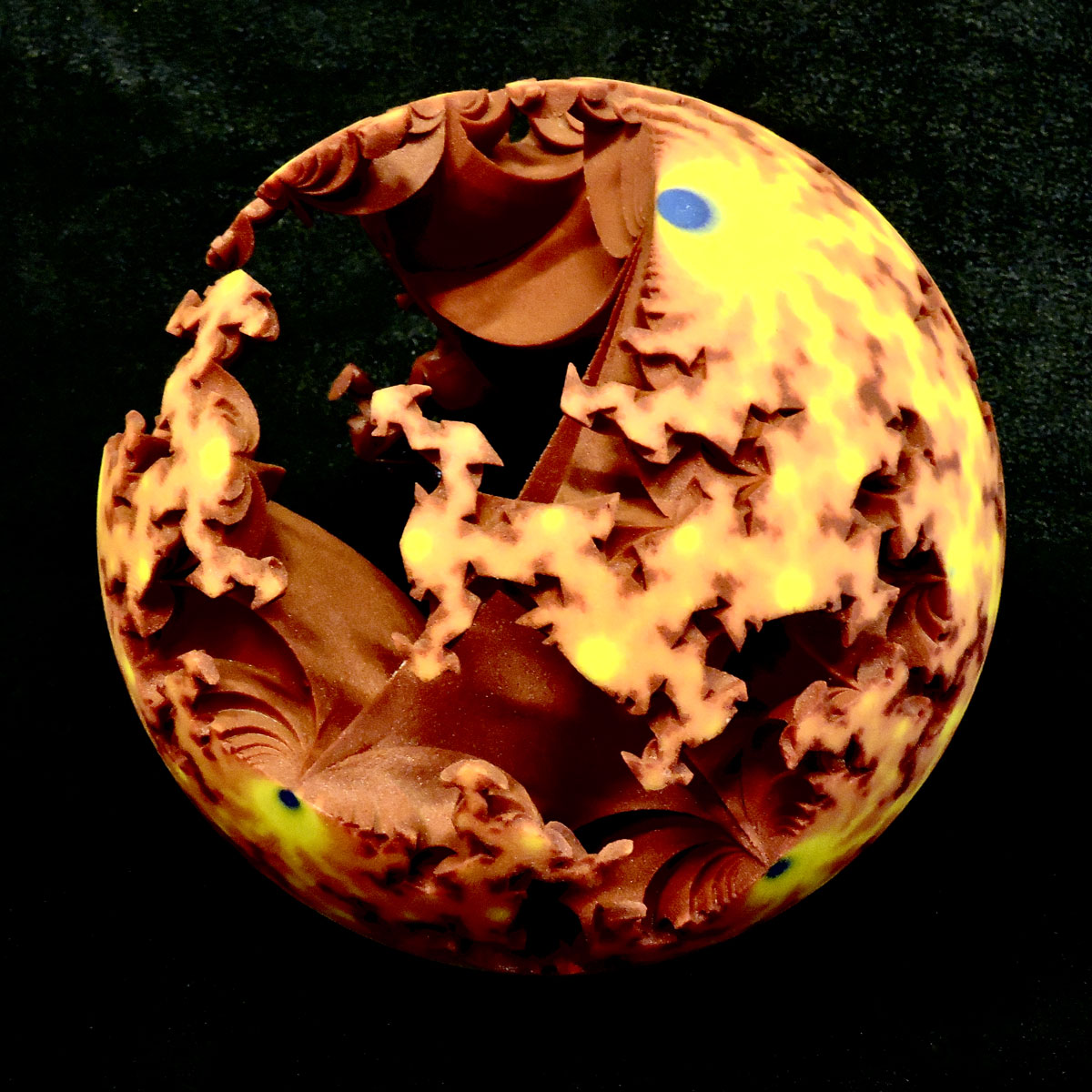}
\caption{\emph{3D printed cohomology fractal and elevation.}}
\label{Fig:3Dprint}
\end{wrapfigure}

\vspace{-20pt}
\hspace{2.2in}
\section*{3D printed sculpture}

In \reffig{3Dprint} we see another depiction of an elevation of the surface as well as the cohomology fractal associated with the figure-eight knot.  This 3D print lies in the Poincar\'e ball model of hyperbolic space.  Half of the ball has been cut away; inside you can see the pleated elevation of the Seifert surface of the knot.  On the remaining boundary we have drawn the cohomology fractal.  Here the fractal is coloured according to a temperature gradient.

\section*{Do these pictures exist? }

That is, as the radius of the visual sphere tends to infinity, is there a ``limiting'' picture?  For example, when drawing many fractals (such as Julia sets), increasing the resolution produces images that converge; this suggests the existence of some underlying mathematical object.  However, as shown in \reffig{VisSphereRadii}, our images do not so converge as the radius tends to infinity.  Instead the limiting object should be a \emph{distribution}.  This cannot be evaluated at a point; rather it is integrated against a characteristic function of a set (for example, a pixel).  We plan to address this in a future paper.

\section*{Acknowledgements}

This work began with our attempts to draw images of Cannon--Thurston maps using pixel techniques and convex geometry.  We thank Fran\c{c}ois Gu\'eritaud for suggesting that we use ray-casting instead.  This material is based in part upon work supported by the National Science Foundation under Grant No. DMS-1439786 and the Alfred P.~Sloan Foundation award G-2019-11406 while the authors were in residence at the Institute for Computational and Experimental Research in Mathematics in Providence, RI, during the Illustrating Mathematics program.  The third author was supported in part by National Science Foundation grant DMS-1708239.  Some work was completed during a research stay at the University of Warwick; this was partially supported by the Simons Foundation and by the Mathematisches Forschungsinstitut Oberwolfach.

\bibliographystyle{hyperplain}
\bibliography{cohomology_fractals}
\end{document}